%% file: main.tex
\documentclass{scrartcl}
\usepackage{graphicx} 
\usepackage{todonotes}
\input{commands}

\title{Finite Simple Groups in the Primitive Positive Constructability Poset}
\author{Sebastian Meyer and Florian Starke\thanks{
Both authors have received funding from the European Research Council (Project POCOCOP, ERC Synergy Grant
101071674). Views and opinions expressed are however those of the authors only and do not necessarily reflect those of the European Union or the European Research Council Executive Agency. Neither the European Union nor the granting authority can be held responsible for them.}}

\date{\today}

\begin{document}

\maketitle

\input{0_abstract}

\input{1_introduction}
\input{2_preliminaries}
\input{3_submaximal}

\input{4_Groups-or-cycles-paper}
\appendix
\input{A_Simple-group-actions}
\newpage
\tableofcontents
\bibliography{citations} 
\bibliographystyle{alpha}
\end{document}

%% file: commands.tex
\usepackage{amsmath}
\usepackage{amsthm}
\usepackage{amssymb}
\usepackage{amsfonts}
\usepackage{mathrsfs}
\usepackage{color}
\usepackage[utf8]{inputenc}
\usepackage{enumitem}
\usepackage{multicol}
\usepackage[UKenglish]{babel}
\usepackage[babel]{csquotes}
\usepackage{mathtools}
\usepackage{rotating}
\usepackage{latexsym}
\usepackage{epsfig}
\usepackage{verbatim}
\usepackage{thmtools,thm-restate}
\usepackage{listings} 

\usepackage{makecell,booktabs,amsxtra}

\usepackage{adjustbox}
\usepackage{ifthen}
\usepackage[position=top]{subfig}
\usepackage{pgfplots}
\pgfplotsset{xticklabels={}, yticklabels={} width=7cm,compat=1.8}
\usepackage{pgfplotstable}
\pgfmathsetseed{1138} 
\pgfplotstableset{ 
	create on use/x/.style={create col/expr={42+2*\pgfplotstablerow}},
	create on use/y/.style={create col/expr={(0.6*\thisrow{x}+130)+5*rand}}
}
\pgfplotstablenew[columns={x,y}]{30}\loadedtable

\usepackage[ruled]{algorithm2e}

\SetAlFnt{\small}
\SetAlCapFnt{\small}
\SetAlCapNameFnt{\small}
\SetAlCapHSkip{0pt}
\IncMargin{-\parindent}

\usepackage{url}
\usepackage[square,numbers,sort]{natbib}

\usepackage{pgfplots,pgfplotstable}
\pgfplotsset{compat=1.8}

\usepackage[bookmarks=false,colorlinks=true,linkcolor=blue,urlcolor=blue,citecolor=red,pdfborder={0 0 0.5}]{hyperref}

\theoremstyle{plain}
\newtheorem{theorem}{Theorem}[section]
\newtheorem{lemma}[theorem]{Lemma}

\newtheorem{corollary}[theorem]{Corollary}

\declaretheoremstyle[
notefont=\bfseries, notebraces={}{},
bodyfont=\normalfont\itshape,
headformat=\NAME \NOTE
]{nopar}
\declaretheorem[style=nopar,name=Corollary]{corollary*}
\declaretheorem[style=nopar,name=Theorem]{theorem*}

\theoremstyle{definition}
\newtheorem{definition}[theorem]{Definition}

\newtheorem*{remark*}{Remark}
\newtheorem{observation}[theorem]{Observation}

\newcount \mycount
\newcount \mycountb

\DeclareMathOperator{\N}{\mathbb N}
\DeclareMathOperator{\Z}{\mathbb Z}

\DeclareMathOperator{\A}{{\mathbb A}}
\DeclareMathOperator{\B}{{\mathbb B}}

\DeclareMathOperator{\bP}{{\mathbb P}}
\DeclareMathOperator{\C}{{\mathbb C}}

\DeclareMathOperator{\T}{{\mathbb T}}
\DeclareMathOperator{\bT}{{\T}}
\renewcommand{\S}{{\mathbb S}}

\renewcommand{\P}{\mathbb P}

\DeclareMathOperator{\Pol}{Pol}

\DeclareMathOperator{\maj}{{\operatorname{maj}}}

\DeclareMathOperator{\Maltsev}{{\operatorname{Mal}}}

\DeclareMathOperator{\Majority}{{\operatorname{Maj}}}

\newcommand{\TS}[1]{{\operatorname{TS}(#1)}}
\newcommand{\FS}[1]{{\operatorname{FS}(#1)}}

\newcommand{\GPinter}[2]{{\operatorname{GP}(#1,#2)}}
\newcommand{\SymGPinter}[1]{{\operatorname{FSGP}(#1)}}

\providecommand{\dotdiv}{
  \mathbin{
    \vphantom{+}
    \text{
      \mathsurround=0pt 
      \protect\ooalign{
        \noalign{\kern-.45ex}
        \hidewidth$\smash{\cdot}$\hidewidth\cr 
        \noalign{\kern.45ex}
        $-$\cr 
      }%
    }%
  }%
}
\providecommand{\cupdot}{
  \mathbin{
    \vphantom{+}
    \text{
      \mathsurround=0pt 
      \protect\ooalign{
        \noalign{\kern-.4ex}
        \hidewidth$\smash{\cdot}$\hidewidth\cr 
        \noalign{\kern.4ex}
        $\cup$\cr 
      }%
    }%
  }%
}

\newcommand{\ppeq}{\mathrel{\equiv}}
\newcommand{\ppleq}{\mathrel{\leq}}

\theoremstyle{definition}
\newtheorem{examplex}[theorem]{Example}
\newenvironment{example}
{\pushQED{\qed}\examplex}
{\popQED\endexamplex}

\theoremstyle{definition}

\makeatletter
\newcommand\footnoteref[1]{\protected@xdef\@thefnmark{\ref{#1}}\@footnotemark}
\makeatother

\makeatletter
\newcommand*{\rom}[1]{\expandafter\@slowromancap\romannumeral #1@}
\makeatother

\makeatletter
\def\blfootnote{\xdef\@thefnmark{}\@footnotetext}
\makeatother

\usepackage{setspace}

\usepackage{epigraph}
\usepackage{siunitx}
\input{inlinetikz}

\newcommand{\pairing}{generalized pairing}

\newcommand{\clone}{\mathcal C}
\newcommand{\cloneD}{\mathcal D}

\newcommand{\Majocolumns}[1]{\maj\begin{pmatrix} #1 \phantom{,}\end{pmatrix}}
\newcommand{\Malcolumns}[1]{\malt\begin{pmatrix} #1 \phantom{,}\end{pmatrix}}
\newcommand{\Malt}{\Maltsev}
\newcommand{\Sym}{\operatorname{Sym}}
\newcommand{\Gmin}{\operatorname{Gmin}} 
\newcommand{\gmin}{g}
\newcommand{\act}{\curvearrowright} 
\newcommand{\prim}{\mathbb{P}}   
\newcommand{\iso}{\cong}    
\newcommand{\IZ}{\mathbb{Z}}    
\newcommand{\IR}{\mathbb{R}}    
\newcommand{\IN}{\mathbb{N}}    
\newcommand{\stab}{\operatorname{Stab}} 
\newcommand{\structA}{\mathbb{A}}
\newcommand{\structB}{\mathbb{B}}

\newcommand{\pol}{\operatorname{Pol}}

\newcommand{\normalsubgroup}{\trianglelefteq}   
\newcommand{\conditionset}{E}
\newcommand{\leftfrac}[2]{#2  \backslash #1}

\DeclareMathOperator{\malt}{{\operatorname{mal}}}
\newcommand{\gp}[2]{p_{#1,#2}}
\newcommand{\fs}[1]{f}
\newcommand{\ts}[1]{t}
\newcommand{\setdifference}{\mathrel{\triangle}}

%% file: inlinetikz.tex
\newcounter{anzahl}
\def\nodeDist{0.4}

\newcommand{\tikzpathFreeVariables}[1]{
\begin{tikzpicture}[scale=0.5]
\setcounter{anzahl}{0}

\foreach \y/\c [count=\xi from 1] in #1 {
    \ifthenelse{\c=0}{
    \node[var-f] (\xi) at (\nodeDist*\xi,\nodeDist*\y) {};}{
    \node[var-b] (\xi) at (\nodeDist*\xi,\nodeDist*\y) {};}
    
    \setcounter{anzahl}{\xi}
    }
    
\foreach \y [count=\yi from 1,count=\yii from 2] in {2,...,\arabic{anzahl}} 
    \path (\yi) edge (\yii);
\end{tikzpicture}
}

\tikzstyle{var-b}=[circle,fill,draw=white,inner sep=0pt,minimum size=3.5pt]
\tikzstyle{bullet}=[circle,fill,draw=white,inner sep=0pt,minimum size=3.2pt]
\tikzstyle{var-f}=[circle,draw,inner sep=0pt,minimum size=3pt]

\usetikzlibrary{arrows,arrows.meta}
\usetikzlibrary{bending,calc}

\newcommand{\Twentyone}{
\begin{tikzpicture}
\node at (0,0) {\TwentyoneC};
\node at (4,0) {\TwentyoneB};
\node at (8,0) {\TwentyoneA};
\end{tikzpicture}
}
\def\scale{0.7}
\newcommand{\TwentyoneA}{
\begin{tikzpicture}[scale=\scale]
\node[bullet, label={right:$17$}] (0) at (2,2) {};
\node[bullet, label={right:$19$}] (1) at (4,3) {};
\node[bullet, label={right:$21$}] (2) at (5,4) {};
\node[bullet, label={right:$16$}] (3) at (2,4) {};
\node[bullet, label={ left:$18$}] (4) at (3,3) {};
\node[bullet, label={right:$20$}] (5) at (5,2) {};
\path[->,>=stealth']
(5) edge (2)
(2) edge (1)
(0) edge (4)
(4) edge (3)
(3) edge (0)
(1) edge (5)
;
\path[dashed]
(5) edge[bend right] (2)
(0) edge[loop,out=135+180,in=45+180,looseness=20] (0)
(4) edge[] (1)
(3) edge[loop,looseness=20] (3)
;
\end{tikzpicture}
}
\newcommand{\TwentyoneB}{
\begin{tikzpicture}[scale=\scale]
\node[bullet, label={right:$7$}] (0) at (1,0) {};
\node[bullet, label={ left:$11$}] (1) at (3,1) {};
\node[bullet, label={ above:$6$}] (2) at (0,0) {};
\node[bullet, label={right:$12$}] (3) at (4,1) {};
\node[bullet, label={south:$15$}] (4) at (3,-2) {};
\node[bullet, label={right:$13$}] (5) at (4,-1) {};
\node[bullet, label={north:$10$}] (6) at (3,2) {};
\node[bullet, label={ left:$8$}] (7) at (2,-1) {};
\node[bullet, label={ left:$9$}] (8) at (2,1) {};
\node[bullet, label={ left:$14$}] (9) at (3,-1) {};
\path[->,>=stealth']
(8) edge (0)
(0) edge (7)
(6) edge (1)
(1) edge (3)
(5) edge (9)
(9) edge (4)
(3) edge (6)
(7) edge (8)
(2) edge[loop,out=135+90,in=45+90,looseness=20] (2)
(4) edge (5)
;
\path[dashed]
(8) edge[] (6)
(5) edge[] (3)
(1) edge[loop,out=135+180,in=45+180,looseness=20] (1)
(0) edge[] (2)
(4) edge[] (7)
(9) edge[loop,looseness=20] (9)
;
\end{tikzpicture}
}

\newcommand{\TwentyoneC}{
\begin{tikzpicture}[scale=\scale]
\node[bullet, label={ left:$1$}] (0) at (0,2) {};
\node[bullet, label={ left:$2$}] (1) at (0,0) {};
\node[bullet, label={right:$3$}] (2) at (1,2) {};
\node[bullet, label={right:$4$}] (3) at (1,0) {};
\node[bullet, label={above:$5$}] (4) at (2,1) {};
\path[->,>=stealth']
(0) edge[loop,out=135+180,in=45+180,looseness=20] (0)
(2) edge (3)
(3) edge (4)
(1) edge[loop,looseness=20] (1)
(4) edge (2)
;
\path[dashed]
(0) edge[] (2)
(1) edge[] (3)
(4) edge[loop,out=135-90,in=45-90,looseness=20] (4)
;
\end{tikzpicture}
}
\newcommand{\TwentytwoA}{
\begin{tikzpicture}
\node[bullet, label={ left:$8$}] (0) at (0,4) {};
\node[bullet, label={ left:$9$}] (1) at (0,3) {};
\node[bullet, label={ left:$10$}] (2) at (-1,2) {};
\node[bullet, label={south:$11$}] (3) at (-2,1) {};
\node[bullet, label={south:$12$}] (4) at (-1,1) {};
\node[bullet, label={south:$13$}] (5) at (1,1) {};
\node[bullet, label={south:$14$}] (6) at (2,1) {};
\node[bullet, label={right:$15$}] (7) at (1,2) {};
\path[->,>=stealth']
(0) edge[loop,out=135+0,in=45+0,looseness=20] (0)
(1) edge (2)
(2) edge (3)
(3) edge (4)
(4) edge (5)
(5) edge (6)
(6) edge (7)
(7) edge (1)
;
\path[dashed]
(0) edge[] (1)
(2) edge[] (7)
(3) edge[bend left] (4)
(5) edge[bend left] (6)
;
\end{tikzpicture}    
}

\newcommand{\TwentytwoB}{
\begin{tikzpicture}[xscale=-1]
\node[bullet, label={south:$1$}] (1) at (0,0) {};
\node[bullet, label={south:$2$}] (2) at (1,0) {};
\node[bullet, label={south:$3$}] (3) at (2,0) {};
\node[bullet, label={south:$4$}] (4) at (3,0) {};
\node[bullet, label={ left:$5$}] (5) at (2,1) {};
\node[bullet, label={north:$6$}] (6) at (1,1) {};
\node[bullet, label={north:$7$}] (7) at (0,1) {};
\path[->,>=stealth']
(1) edge (2)
(2) edge (3)
(3) edge (4)
(4) edge (5)
(5) edge (6)
(6) edge (7)
(7) edge (1)
;
\path[dashed]
(1) edge[loop,out=135+90,in=45+90,looseness=20] (1)
(5) edge[loop,out=135+0,in=45+0,looseness=20] (5)
(7) edge[loop,out=135+90,in=45+90,looseness=20] (7)
(2) edge[] (6)
(3) edge[bend left] (4)
;

\node[bullet,white] (0) at (0,3) {};
\path[->,>=stealth',white]
(0) edge[loop,looseness=20] (0)
;
\end{tikzpicture}    
}
\newcommand{\TwentytwoC}{
\begin{tikzpicture}[xscale=-1]
\node[bullet, label={south:$22$}] (1) at (-0,0) {};
\node[bullet, label={south:$21$}] (2) at (-1,0) {};
\node[bullet, label={south:$20$}] (3) at (-2,0) {};
\node[bullet, label={south:$19$}] (4) at (-3,0) {};
\node[bullet, label={right:$18$}] (5) at (-2,1) {};
\node[bullet, label={north:$17$}] (6) at (-1,1) {};
\node[bullet, label={north:$16$}] (7) at (-0,1) {};
\path[->,>=stealth']
(1) edge (7)
(2) edge (1)
(3) edge (2)
(4) edge (3)
(5) edge (4)
(6) edge (5)
(7) edge (6)
;
\path[dashed]
(1) edge[loop,out=135-90,in=45-90,looseness=20] (1)
(5) edge[loop,out=135+0,in=45+0,looseness=20] (5)
(7) edge[loop,out=135-90,in=45-90,looseness=20] (7)
(2) edge[] (6)
(3) edge[bend right] (4)
;
\node[bullet,white] (0) at (0,3) {};
\path[->,>=stealth',white]
(0) edge[loop,looseness=20] (0)
;
\end{tikzpicture}    
}

\newcommand{\FiftytwoA}{
\begin{tikzpicture}
\node[bullet, label={ left:$7$}] (0) at (-1,2) {};
\node[bullet, label={ left:$8$}] (1) at (-1,1) {};
\node[bullet, label={south:$9$}] (2) at (-1,0) {};
\node[bullet, label={south:$10$}] (3) at (0,0) {};
\node[bullet, label={north:$11$}] (4) at (0,1) {};
\node[bullet, label={north:$12$}] (5) at (1,1) {};
\node[bullet, label={south:$13$}] (6) at (1,0) {};
\node[bullet, label={south:$14$}] (7) at (2,0) {};
\node[bullet, label={right:$15$}] (8) at (2,1) {};
\node[bullet, label={right:$16$}] (9) at (2,2) {};
\path[->,>=stealth']
(0) edge[loop,looseness=20] (0)
(1) edge (2)
(2) edge (3)
(3) edge (4)
(4) edge (1)
(5) edge (6)
(6) edge (7)
(7) edge (8)
(8) edge (5)
(9) edge[loop,looseness=20] (9)
;
\path[dashed]
(0) edge[] (1)
(2) edge[loop,out=135+90,in=45+90,looseness=20] (2)
(3) edge[] (5)
(4) edge[] (6)
(7) edge[loop,out=135-90,in=45-90,looseness=20] (7)
(8) edge[] (9)
;
\end{tikzpicture}    
}
\newcommand{\FiftytwoB}{
\begin{tikzpicture}
\node[bullet, label={south:$37$}] (0) at (6,1) {};
\node[bullet, label={right:$33$}] (1) at (6,2) {};
\node[bullet, label={right:$34$}] (2) at (6,3) {};
\node[bullet, label={north:$31$}] (3) at (5,3) {};
\node[bullet, label={south:$32$}] (4) at (5,2) {};
\node[bullet, label={south:$28$}] (5) at (4,3) {};
\node[bullet, label={right:$29$}] (6) at (4,4) {};
\node[bullet, label={ left:$30$}] (7) at (3,4) {};
\node[bullet, label={south:$27$}] (8) at (3,3) {};
\node[bullet, label={north:$26$}] (9) at (2,3) {};
\node[bullet, label={ left:$23$}] (10) at (1,3) {};
\node[bullet, label={ left:$24$}] (11) at (1,2) {};
\node[bullet, label={south:$25$}] (12) at (2,2) {};
\node[bullet, label={south:$35$}] (13) at (3,1) {};
\node[bullet, label={south:$36$}] (14) at (4,1) {};
\path[->,>=stealth']
(0) edge[loop,out=135+90,in=45+90,looseness=20] (0)
(1) edge (2)
(2) edge (3)
(3) edge (4)
(4) edge (1)
(5) edge (6)
(6) edge (7)
(7) edge (8)
(8) edge (5)
(9) edge (10)
(10) edge (11)
(11) edge (12)
(12) edge (9)
(13) edge (14)
(14) edge (13)
;
\path[dashed]
(0) edge[] (1)
(2) edge[bend right=90] (10)
(3) edge[] (5)
(4) edge[] (14)
(6) edge[loop,looseness=20] (6)
(7) edge[loop,looseness=20] (7)
(8) edge[] (9)
(11) edge[] (13)
(12) edge[loop,out=135-90,in=45-90,looseness=20] (12)
;
\end{tikzpicture}    
}
\newcommand{\FiftytwoBMirrow}{
\begin{tikzpicture}
\node[bullet, label={south:$52$}] (0) at (-6,1) {};
\node[bullet, label={ left:$39$}] (1) at (-6,2) {};
\node[bullet, label={ left:$38$}] (2) at (-6,3) {};
\node[bullet, label={north:$41$}] (3) at (-5,3) {};
\node[bullet, label={south:$40$}] (4) at (-5,2) {};
\node[bullet, label={south:$42$}] (5) at (-4,3) {};
\node[bullet, label={ left:$45$}] (6) at (-4,4) {};
\node[bullet, label={right:$44$}] (7) at (-3,4) {};
\node[bullet, label={south:$43$}] (8) at (-3,3) {};
\node[bullet, label={north:$46$}] (9) at (-2,3) {};
\node[bullet, label={right:$49$}] (10) at (-1,3) {};
\node[bullet, label={right:$48$}] (11) at (-1,2) {};
\node[bullet, label={south:$47$}] (12) at (-2,2) {};
\node[bullet, label={south:$51$}] (13) at (-3,1) {};
\node[bullet, label={south:$50$}] (14) at (-4,1) {};
\path[<-,>=stealth']
(0) edge[loop,out=135-90,in=45-90,looseness=20] (0)
(1) edge (2)
(2) edge (3)
(3) edge (4)
(4) edge (1)
(5) edge (6)
(6) edge (7)
(7) edge (8)
(8) edge (5)
(9) edge (10)
(10) edge (11)
(11) edge (12)
(12) edge (9)
(13) edge (14)
(14) edge (13)
;
\path[dashed]
(0) edge[] (1)
(2) edge[ bend left=90] (10)
(3) edge[] (5)
(4) edge[] (14)
(6) edge[loop,looseness=20] (6)
(7) edge[loop,looseness=20] (7)
(8) edge[] (9)
(11) edge[] (13)
(12) edge[loop,out=135+90,in=45+90,looseness=20] (12)
;
\end{tikzpicture}    
}

\newcommand{\FiftytwoC}{
\begin{tikzpicture}
\node[bullet, label={south:$2$}] (1) at (0,0) {};
\node[bullet, label={south:$3$}] (2) at (1,0) {};
\node[bullet, label={right:$5$}] (3) at (1,1) {};
\node[bullet, label={ left:$1$}] (4) at (0,1) {};
\node[bullet, label={south:$4$}] (5) at (2,0) {};
\node[bullet, label={right:$6$}] (6) at (2,1) {};
\path[->,>=stealth']
(1) edge (2)
(2) edge (3)
(3) edge (4)
(4) edge (1)
(5) edge[<->] (6)
;
\path[dashed]
(1) edge[ bend right] (4)
(2) edge[] (5)
(6) edge[loop,looseness=20] (6)
(3) edge[loop,looseness=20] (3)
;
\end{tikzpicture}    
}

\newcommand{\FiftytwoCMirrow}{
\begin{tikzpicture}
\node[bullet, label={south:$20$}] (1) at (-0,0) {};
\node[bullet, label={south:$19$}] (2) at (-1,0) {};
\node[bullet, label={ left:$22$}] (3) at (-1,1) {};
\node[bullet, label={right:$21$}] (4) at (-0,1) {};
\node[bullet, label={south:$18$}] (5) at (-2,0) {};
\node[bullet, label={ left:$17$}] (6) at (-2,1) {};
\path[<-,>=stealth']
(1) edge (2)
(2) edge (3)
(3) edge (4)
(4) edge (1)
(5) edge[<->] (6)
;
\path[dashed]
(1) edge[bend left] (4)
(2) edge[] (5)
(6) edge[loop,looseness=20] (6)
(3) edge[loop,looseness=20] (3)
;
\end{tikzpicture}    
}

%% file: 0_abstract.tex
\abstract{
Primitive positive constructions between structures are motivated by computational complexity of constraint satisfaction problems. They induce a partial order on equivalence classes of finite structures. We show that this order has an infinite third layer, given by the equivalence class of an oriented graph and of infinitely many permutation groups that are abstractly isomorphic to all the finite simple groups in a one-to-one correspondence.
More concretely, the two topmost layers of the primitive positive constructability poset are known to consist of a single element each. We consider $\P_1$, a representative of the element in the second layer, and its polymorphism clone $\Pol(\P_1)$.
We show that any clone over a finite domain that admits a quasi Maltsev operation and fully symmetric operations of all arities admits an incoming minion homomorphism from $\Pol(\P_1)$ and use this result to show that in the primitive positive contructability poset, the lower covers of $\P_1$ are represented by the transitive tournament on three vertices and for each finite simple group by the disjoint union of all its primitive group actions. 
}

%% file: 1_introduction.tex
\section{Introduction}



The primitive positive constructability poset is closely linked to the classification of clones in universal algebra. 
A \emph{clone} is a set of operations (over a fixed domain) that contains all projections and is closed under composition. The set of all clones over a fixed domain, ordered by inclusion, forms a lattice. In 1941 Post gave a description of the lattice of all clones on a two-element domain~\cite{Post}. In 1959 Janov and Mučnik showed that there exist  continuum many clones over a $k$-element domain, for every $k \geq 3$~\cite{3elem}. The arising lattices are considerably more complex and a full description of them as in the two-element case is unlikely. 
Note that each of these lattices has the clone of all operations as largest element and the clone containing only projections as smallest element.  Remarkably, Rosenberg managed to show in 1965 that the largest clone always has six lower covers (if the domain has at least three elements)~\cite{RosenbergMaximal}. 

Since the late 90's clones have gained traction in the area of \emph{constraint satisfaction problems} (CSPs) \cite{Jeavons}: Let $\A$ be a finite relational structure, the CSP of $\A$ is the following computational problem.
\begin{align*}
    &\textbf{Input:} \text{ a finite structure $\B$ (with the same signature as $\A$)}\\
    &\textbf{Output:}\text{ \textbf{Yes} if $\B$ has a homomorphism to $\A$, \textbf{No} otherwise}
\end{align*}
It turns out that the computational complexity (up to logspace reductions) of the CSP of $\A$  only depends on the minor conditions satisfied by the polymorphisms of $\A$~\cite{wonderland}. Roughly speaking, the more minor conditions are satisfied by the polymorphisms of $\A$, the easier the CSP is. For example, if $\A$ has a ternary majority polymorphism, then the CSP of $\A$ is solved by a path consistency algorithm in polynomial time \cite{FederVardi}.
The polymorphisms of $\A$ form a clone over the domain $A$. 
Therefore, it is interesting to consider clones ordered by the minor conditions that they satisfy: 
For two structures $\A$ and $\B$ with a finite domain, we define $\A\ppleq \B$ if 
all minor conditions satisfied by $\Pol(\A)$ are also satisfied by $\Pol(\B)$. This is the case if and only if there is a \emph{minion homomorphism} from $\pol(\A)$ to $ \pol(\B)$, or, equivalently, if $\A$ \emph{primitively positively (pp-) constructs} $\B$ (see Theorem~\ref{thm:ppconIffMinionHom}).
This pre-order defines a poset on the quotient of all structures with finite domains up to pp-interconstructability. We call this poset the \emph{pp-constructability poset} (of all finite structures).  
This poset has been studied before~\cite{wonderland,vucajBodirskytwoElementPPPoset,smooth-digraphs,maximal-digraphs,vucaj2024submaximal}. 
It is well known that all clones over a finite domain are polymorphism clones of a suitable finite structure. Therefore, we can define the same poset with clones over a finite domain instead of structures with a finite domain. Note that the corresponding pre-order on clones extends the inclusion order and that we can now compare clones over different domains. Hence, there is hope that studying the pp-constructability poset also leads to a better understanding of clones on finite domains. 

The pp-constructability poset is known to be infinite and  of size at most $|\IR|=2^{\aleph_0}$. It has a greatest element, represented by $\C_1$, a single vertex with a loop. It has a second layer, consisting of a single element which can be represented by $\P_1$, the digraph consisting of a single directed edge. Moreover, the poset has a smallest element given by the equivalence class of the complete graph on three vertices. 
It is known that there is no second-to-last layer but infinite descending chains approaching the minimal element. 
Studying this poset for all finite structures is very difficult. However, the induced subposet consisting of all structures with at most two elements and the subposet consisting of all smooth digraphs are completely known and are presented in \cite{vucajBodirskytwoElementPPPoset} and \cite{smooth-digraphs}, respectively.
This article is split in two parts. In the first part we study $\mathcal I$, the clone of all idempotent operations on a two-element domain. This clone is 
the polymorphism clone of $\P_1$. In particular, we are interested in a generating set of the set of minor conditions that this clone satisfies, i.e., a set of minor conditions $M$ such that any clone over a finite domain that satisfies all conditions in $M$ also satisfies all conditions satisfied by $\mathcal{I}$. Vucaj and Zhuk showed that the set consisting of the totally symmetric condition for each arity and the generalized minority condition for each odd arity forms such a generating set~\cite{vucaj2024submaximal}. They conjectured that the set would still be generating if the generalized minority conditions would be replaced by the single quasi Maltsev condition. In Theorem~\ref{thm:idempotence} we show that this conjecture is true and that a set of even weaker conditions is also a generator, namely the set consisting of the quasi Maltsev condition and of the fully symmetric conditions of all arities. Later in Corollaries~\ref{cor:weakerConditionBasisForP1} and~\ref{cor:weakestGeneratingSet} we show that the fully symmetric conditions can be replaced by slightly weaker conditions such that the resulting generating set is, in a certain sense, minimal.

In the second part of this article we use the results from the first part to study the third layer of the poset of all finite structures ordered by pp-constructability. 
For digraphs ordered by pp-constructability Bodirsky and Starke showed that the third layer of the poset consists of the digraphs $\T_3$, i.e., the transitive tournament on three vertices, and all directed cycles of prime length~\cite{maximal-digraphs}. The order on these digraphs is as depicted below
\begin{center}
        \begin{tikzpicture}[scale=1.0]
    \node (0) at (2,2)  {$\C_{1}$};
    \node (1) at (2,1)  {$\P_{1}$};
    \node (20) at (0,0)  {$\T_{3}$};
    \node (21) at (1,0)  {$\C_{2}$};
    \node (22) at (2,0)  {$\C_{3}$};
    \node (23) at (3,0)  {$\C_{5}$};
    \node (24) at (4,0)  {$\dots$};
    \path
        (0) edge (1)
        (1) edge (20)
        (1) edge (21)
        (1) edge (22)
        (1) edge (23)
        (1) edge (24)
        ;
    \end{tikzpicture}
\end{center}

Recently, Bodirsky and Starke showed that $\T_3$ is also a lower cover of $\P_1$ in the poset of all finite structures~\cite{bodirskyStarke2024symmetriclineararcmonadicUnfoldedCaterpillar}. 
We will use the results from the first part of this article about the polymorphism clone of $\P_1$ to show that the directed cycles of prime lengths are lower covers of $\P_1$ in the poset of all finite structures and that they correspond to specific group actions of abelian finite simple groups. 
We will show that there are other lower covers of $\P_1$ and that they are the structures corresponding to specific actions of non-abelian finite simple groups. Surprisingly, one of this structures has been studied before in a different context. In  \cite{carvalhoKrokhin21ElementStructure} the authors describe a 21-element structure which they found by the help of a computer. They were interested in this structure because it was the first known example of a structure that has cyclic polymorphisms of all arities and does not have fully symmetric polymorphisms of all arities. It turns out that the structure they found is closely linked to the smallest non-abelian finite simple group $A_5$ and one of the structures on the third layer.
Using the notation from Definition~\ref{def:S(GactX)},
we show that the first three layers of the pp-constructability poset on all finite structures look like this:
\begin{center}
    \centering
        \begin{tikzpicture}[scale=1.0]
    \node (0) at (1,2)  {$\C_{1}$};
    \node (1) at (1,1)  {$\P_{1}$};
    \node (20) at (0,0)  {$\T_{3}$};
    \node (21) at (1,0)  {\phantom{A} };
    \node (22) at (2,0)  {\phantom{A}};
    \node (23) at (3,0)  {\phantom{AAAABBXCCC} };
    \node (24) at (4,0)  {$\{\S(G\act\prim(G))\mid G \text{ finite simple group}\}$};
    \path
        (0) edge (1)
        (1) edge (20)
        (1) edge (21)
        (1) edge (22)
        (1) edge (23)
        (1) edge (24)
        ;
    \end{tikzpicture}
\end{center}
This result can be found in Theorem~\ref{thm:mainLowerCoversClassification}. In particular, our result implies that the classification of finite simple groups is relevant when studying the pp-constructability poset. We use the classification theorem only for a comment in the final section.


Our result also has applications in descriptive complexity theory. 
Recently, Bodirsky and Starke presented symmetric linear arc monadic Datalog, a natural Datalog fragment that solves exactly the CSPs of structures that can be pp-constructed from $\P_1$ \cite{bodirskyStarke2024symmetriclineararcmonadicUnfoldedCaterpillar}. Our results show that the structures whose CSP cannot be solved by this Datalog fragment are exactly those that can pp-construct either $\T_3$ or $\S(G\act\prim(G))$ for some finite simple group $G$. 

Moreover, it is known that a CSP of a structure $\structA$ can be solved by basic linear programming robustly in polynomial time if and only if $\structA$ has fully symmetric polymorphisms of all arities \cite{BLP}. It follows from our results that this is equivalent to $\structA$ not pp-constructing $\S(G\act\prim(G))$ for any finite simple group $G$. 

%% file: 2_preliminaries.tex
\section{Preliminaries}
We write $[n]$ for the set $\{1,\dots,n\}$.


\subsection{Group actions}
In this section we review some basic notation for groups and group actions. 
For a set $X$ we denote the \emph{symmetric group of $X$} by $\Sym(X)$. For $n\in\N$ we let $S_n$ denote $\Sym([n])$. 
A \emph{group action} of a group $G$ on a set $X$ is a homomorphism $\alpha\colon G\to\Sym(X)$. This action is denoted by $G\curvearrowright_\alpha X$ or $G\curvearrowright X$ if $\alpha$ is clear from the context. For $g\in G$ and $x\in X$ we often write $gx$ or $g.x$ instead of $(\alpha(g))(x)$. 
A \emph{$G$-set} is a set together with a group action of $G$ on it.
Common $G$-sets that will be relevant in this article are
\begin{itemize}
    \item $G$ itself, where the action is given by left multiplication,
    \item the left cosets $G/H=\{gH\mid g\in G\}$ with respect to a subgroup $H$, where the action is given by left multiplication,
    \item the set of maximal subgroups $\{M< G\mid M\text{ maximal subgroup}\}$, where the action is given by conjugation as $g.M=gMg^{-1}$ for $g\in G$,
    \item $X$ if $G\leq\Sym(X)$, in particular $[n]$ as $S_n$-set
    \item $X^Y$, the set of all functions from $Y$ to $X$, if $X$ is a $G$-set. The action is given by post-composition (or componentwise action), i.e., 
    $$g.t=(Y\to X,y\mapsto g(t(y))).$$
    \item $X^Y$, the set of all functions from $Y$ to $X$,
    if $Y$ is a $G$-set. The action is given by pre-composition, i.e., 
    $$g.t=t_g=(Y\to X,y\mapsto t(gy)).$$ 
    Note that this is a right action, i.e., $(gh).t=h.(g.t)$ respectively $t_{gh}=(t_g)_h$.
\end{itemize}


\begin{definition} \label{def:primitiveAction}
    A group action $G \act X$ is called \emph{primitive} if $X$ has more than one element and the only  partitions of $X$ that are respected by the $G$ action are $\{X\}$ (a single partition) and $\{\{x\}\mid x \in X\}$ (the discrete partition). A $G$-set is \emph{primitive} if its action is.
\end{definition}

Let $G\act X$ be a group action and $S\subseteq X$. Then the \emph{(set) stabiliser of $S$ with respect to $G$} is denoted by $\stab_G(S)$. 
The \emph{orbit of $S$ under $G$} is denoted by $G(S)$. 
If $G$ is clear from the context we often write $\stab(S)$ instead of $\stab_G(S)$. If $S=\{s\}$ we often write $\stab_G(s)$ instead of $\stab_G(\{s\})$ and $G(s)$ instead of $G(\{s\})$. 

\begin{observation}
    For the group action $G\act G/H$ on the left-$H$-cosets, we have $\stab_G(\{H\})=H$. Conversely, if $G\act X$ is a group action with a single orbit and $x\in X$, then $X$ is isomorphic to $G/\stab_G(x)$ as $G$-set, where an isomorphism is given by $G/\stab_G({x})\to X, g\stab_G({x})\mapsto g.x$. This map is well-defined.

    For different points $x$ and $g.x$ in the same orbit of a $G$-action $G\act X$, we have
    $$
        \stab_G(g.x)=g\stab_G(x)g^{-1},
    $$
    so the stabilizers are conjugated. Similarly, the $G$-sets $G/H$ and $G/H'$ are isomorphic if and only if $H$ and $H'$ are conjugated.
\end{observation}

\begin{observation}\label{obs:primitiveActionsOnMaximalSubgroups}
It is well known that primitive actions are closely related to maximal subgroups. The following are equivalent for a $G$-action $G\act X$:
\begin{enumerate}
    \item $G\act X$ is primitive.
    \item $G\act X$ is isomorphic to $G/M$ as $G$-set for a maximal subgroup $M<G$.
    \item $G\act X$ is transitive and there exists $x\in X$, such that $\stab_G(x)$ is a maximal subgroup of $G$.
    \item $G\act X$ is transitive and for all $x\in X$, $\stab_G(x)$ is a maximal subgroup of $G$.
\end{enumerate}
The $G$-sets $G/M$ and $G/M'$ are isomorphic if and only if $M$ and $M'$ are conjugated. See for example \cite[Corollary 1.4A and 1.5A]{DixonMortimer}.
\end{observation}

\subsubsection{Biactions}
We will end this section with a lemma we will need in the proof of Lemma~\ref{lem:SimpleGroupsModelSimpleConditions}. We consider a set with a right and a left group action. In this case, it is possible to link the two groups using stabilizers.

For a group $G$ with subgroup $H$ we denote the \emph{right cosets of $H$} by $\leftfrac{G}{ H}$. 

\begin{lemma}\label{lem:IsoOfSubquotient}
    Let $G\act X$ and $H\act Y$ be group actions and let $t\in X^Y$. Let $Z_t=G(t)\cap H(t)\subseteq X^Y$. 
    Then, the maps
    \begin{align*}
        \stab_G(H(t))/\stab_G(t) &\to Z_t \text{ and} & \leftfrac{\stab_H(G(t)) }{\stab_H(t)}&\to Z_t\\
        g\stab_G(t) &\mapsto g.t & \stab_H(t) h &\mapsto t_h
    \end{align*}
    are bijections. 
    Moreover, $\stab_G(t)\normalsubgroup \stab_G(H(t))\le G$ and $\stab_H(t)\normalsubgroup \stab_H(G(t))\le H$  
    and the above maps induce an isomorphism
    $$
        \stab_G(H(t)) / \stab_G(t) \iso \stab_H(G(t)) /\stab_H(t)
    $$
    of groups.
\end{lemma}

\begin{proof}
    First we show that $\stab_G(t)\le \stab_G(H(t))$.
    Let $g\in \stab_G(t)$. For $r=t_h\in H(t)$, we have $gr=gt_h=(gt)_h=t_h=r$, thus $g\in \stab_G(H(t))$.

    The map $f\colon\stab_G(H(t))/\stab_G(t) \to Z_t, g\stab_G(t) \mapsto g.t$ is well-defined, since the action of $\stab_G(t)$ on $t$ is by definition trivial. 
    To see why $f$ maps into $Z_t$ consider $g\in \stab_G(H(t))$. By definition there is some $h\in H$ with $gt=t_h$. Hence, 
    $$f(g\stab_G(t))=gt=t_h\in G(t)\cap H(t)=Z_t.$$ 
    The map $f$ is injective, since $gt=g't$ implies $g^{-1}g'\in\stab_G(t)$. 
    It is surjective, since $gt=t_h\in Z_t$ implies that for every $r=t_{h'}\in H(t)$ we have \[gr=g(t_{h'})=(gt)_{h'}=(t_h)_{h'}=t_{hh'}\in H(t)\] and thus $g\in \stab_G(H(t))$.
    
    Now we show that $\stab_G(t)\normalsubgroup \stab_G(H(t))$.
    Let $g'\in \stab_G(t)$ and $g\in \stab_G(H(t))$. Since $f$ maps into $Z_t$ there is an $h\in H$ such that $gt=t_h$. Hence
    $$
        g^{-1}g'gt=g^{-1}g't_h=g^{-1}(g't)_h=g^{-1}t_h=(g^{-1}t)_h=t_{h^{-1}h}=t
    $$
    and $g^{-1}g'g \in \stab_G(t)$ as desired.
    The proof for $H$ with the right action is analogous.

    To see that the quotients are isomorphic, note that the bijections induce two group structures on $Z_t$, which we call $\ast_G$ and $\ast_H$. Note that they have the same unit $1_Gt=t=t_{1_H}$ and $gt \ast_G r=gr$ and $r\ast_Ht_h=r_h$ by definition. Moreover, $gt\ast_Gt_h=g(t_h)=(gt)_h=gt\ast_Ht_h$ for all $gt,t_h\in Z_t$. Hence, both products are equal on $Z_t$. Thus, the groups are isomorphic.
\end{proof}

\subsection{Structures}
We assume familiarity with the concepts of relational structures and first-order formulas from mathematical logic, as introduced for instance in~\cite{Hodges}. 
A \emph{(directed) graph} is a structure with one binary relation. The following graphs will be important in this article, see also Figure~\ref{fig:graphsToStructures}.
\begin{itemize}
    \item The graph $\C_1$ with a single vertex and a loop at this vertex.
    \item The graph $\bP_1$ with two vertices and a single edge between them. 
    \item The transitive tournament $\T_3$ on three vertices $\{0,1,2\}$, which is the graph having the edges $\{(0,1),(0,2),(1,2)\}$.
    \item The cyclic graph $\C_p$ on $p$ vertices $\{0,1,\dots,p-1\}$, which has the edges $\{(0,1),(1,2),\dots,(p-2,p-1),(p-1,0)\}$.
\end{itemize}
In this article, we consider finite structures which are for us structures whose domain is finite. The signature is allowed to be infinite, but it will turn out to be finite in all relevant cases.

\begin{figure}
    \centering
    $$
    \begin{matrix}
    \begin{tikzpicture}[scale=\scale]
        \node[bullet, label={}] (1) at (0,0) {};
        \path[->,>=stealth'] (1) edge[loop,looseness=20] (1);
    \end{tikzpicture} &\hspace{20mm}&
    \begin{tikzpicture}[scale=\scale]
        \node[bullet, label={}] (1) at (0,0) {};
        \node[bullet, label={}] (2) at (1,0) {};
        \path[->,>=stealth'] (1) edge (2);
    \end{tikzpicture} &\hspace{20mm}&
    \begin{tikzpicture}[scale=\scale]
        \node[bullet, label={}] (1) at (0,0) {};
        \node[bullet, label={}] (2) at (1,1) {};
        \node[bullet, label={}] (3) at (0,2) {};
        \path[->,>=stealth'] (1) edge (2) (1) edge (3) (2) edge (3);
    \end{tikzpicture} &\hspace{20mm}&
    \begin{tikzpicture}[scale=\scale]
        \node[bullet, label={}] (1) at (0:1) {};
        \node[bullet, label={}] (2) at (72:1) {};
        \node[bullet, label={}] (3) at (144:1) {};
        \node[bullet, label={}] (4) at (-144:1) {};
        \node[bullet, label={}] (5) at (-72:1) {};
        \path[->,>=stealth'] (1) edge (2) (2) edge (3) (3) edge (4) (4) edge (5) (5) edge (1);
    \end{tikzpicture}
    \\
    \C_1 
    && \P_1 
    && \T_3 
    && \C_5
    \end{matrix}
    $$
    \caption{Some important structures from this article.}
    \label{fig:graphsToStructures}
\end{figure}


\subsection{Clones}
Let $C$ be a set. An \emph{operation on $C$} is a map from $C^n$ to $C$ for some $n\in\N\setminus \{0\}$. 
Let $f\colon C^n\to C$, and $g_1,\dots,g_n\colon C^m\to C$. Then the \emph{composition} of $f$ and $(g_1,\dots,g_n)$, denoted $f(g_1,\dots,g_n)$, is defined as the map
\begin{align*}
    C^m&\to C\\
    (c_1,\dots,c_m)&\mapsto f(g_1(c_1,\dots,c_m),\dots,g_n(c_1,\dots,c_m)).
\end{align*}
For $n\in\N^+$ and $1\leq k\leq n$ the map $\pi_k^n\colon C^n\to C,(c_1,\dots,c_n)\mapsto c_k$ is the \emph{$n$-ary $i$-th projection of $C$}. 
    A \emph{clone} on the set $C$ is a set $\clone$, which can be written as a disjoint union $\coprod_{n\in \IN} \clone^{(n)}$ such that $\clone^{(n)}$ consists of  maps from $C^n$ to $C$, $\clone$ contains all projections of $C$, and $\clone$ is closed under composition.

For $n\in\N^+$ an $n$-ary \emph{polymorphism} of a structure $\structA$ is a homomorphism from $\structA^n$ to $\structA$.
The set of all polymorphisms of $\structA$ is denoted by $\pol(\structA)$. This set is called the \emph{polymorphism clone} of $\structA$. Note that $\pol(\structA)$ contains all projections and is closed under composition. Hence it is indeed a clone.
It is well known that each clone on a finite set is the polymorphism clone of a suitable finite structure \cite{Geiger}.

An operation $f$ on $C$ is called \emph{idempotent}, if $f(c,\dots,c)=c$ for all $c\in C$. A clone is called idempotent if all of its operations are idempotent. Note that all projections are idempotent. 



\subsection{Minor Conditions}\label{sec:minorconditions}
It would be natural to consider the category of clones together with clone homomorphisms. However, it turns out that the wider class of minor preserving maps gives an easier classification while still catching the important differences.

Let $\clone$ be a clone and $f\in\clone^{(n)}$. For any $m\in\N^+$ and any map $\sigma\colon [n]\to [m]$ the map
\begin{align*}
    f_\sigma\colon C^m&\to C, (c_1,\dots,c_m)\mapsto f(c_{\sigma(1)},\dots,c_{\sigma(n)})
\end{align*}
is called a \emph{minor} of $f$. Note that $f_\sigma$ can be obtained from $f$ by composition with projections. Hence, $f_\sigma$ is an element of $\clone^{(m)}$. 
Let $\clone$ and $\cloneD$ be clones. A map $\xi\colon\clone\to\cloneD$ is \emph{minor preserving} (or a \emph{minion homomorphism}) if $\xi(\clone^{(n)})\subseteq\cloneD^{(n)}$ for all $n\in\N^+$ and $\xi(f_\sigma)=\xi(f)_\sigma$ for all $n,m\in\N^+$, all $f\in\clone^{(n)}$, and all maps $\sigma\colon[n]\to[m]$.

\begin{definition}
    A \emph{minor identity} in $(f,g)$ is an identity of the form
    \[
        \forall x_1,\dots,x_n\colon f(y_1,\dots,y_m)=g(z_1,\dots,z_k)
    \]
    where $y_1,\dots,y_m,z_1,\dots,z_k$ are elements of $\{x_1,\dots,x_n\}$ and $f$ and $g$ are $m$ and $k$ ary function symbols, respectively. Note that there are maps $\sigma\colon [m]\to [n]$ and $\tau\colon [k]\to [n]$ such that this minor identity can be written as
    \[
        \forall x_1,\dots,x_n\colon  f_\sigma(x_1,\dots,x_n)=g_\tau(x_1,\dots,x_n)
    \]
    this justifies the name minor identity.  
    A \emph{minor condition} in  $(f_1,\dots,f_n)$ is a  finite set of minor identities in tuples of elements of $\{f_1,\dots,f_n\}$. 
    
    For a minor condition $\Sigma$ in $(f_1,\dots,f_n)$ and a set $\mathcal F$ of operations on the set $C$ we say that $\mathcal F$ \emph{satisfies} $\Sigma$, denoted $\mathcal F \models \Sigma$, if there are $\tilde f_1,\dots,\tilde f_n\in\mathcal{F}$ such that for each minor identity 
    \[
        \forall x_1,\dots,x_n\colon f_i(y_1,\dots,y_m)=f_j(z_1,\dots,z_k)
    \] in $\Sigma$ we have that 
    \[
         \tilde f_i(y_1,\dots,y_m)=\tilde f_j(z_1,\dots,z_k)\text{ for all }x_1,\dots,x_n\in C.
    \]
    We write $f\models \Sigma$ for $\{f\}\models \Sigma$.
    Let $\Sigma_1$ and $\Sigma_2$ be two minor conditions. We say that $\Sigma_1$ \emph{implies $\Sigma_2$ (for clones)} if any clone $\clone$ that satisfies $\Sigma_1$ also satisfies $\Sigma_2$. Two minor conditions are \emph{equivalent (for clones)} if they imply each other.
\end{definition}

\begin{definition}
We consider these minor conditions.
\begin{itemize}
\item The \emph{quasi Maltsev} condition $\Malt$ consists of the identities
    $$
        \forall x,y \colon f(x,x,x)=f(x,y,y)=f(y,y,x)
    $$
    \item The \emph{quasi majority} condition $\Majority$ consists of the identities
    $$
        \forall x,y \colon f(x,x,x)=f(x,x,y)=f(x,y,x)=f(y,x,x)
    $$
    \item The \emph{$p$-cyclic} condition $\Sigma_p$ consists of the identity

    $$
        \forall x_1,\dots,x_p\colon f(x_1,\dots,x_{p-1},x_p)=f(x_2,\dots,x_p,x_1)
    $$
    
    \item The \emph{fully symmetric} condition $ \FS n $ contains all identities 
    $$
        \forall x_1,\dots,x_n\colon f(x_1,\dots,x_n)=f(y_1,\dots,y_n)
    $$
    where $x_1,\dots,x_n,y_1,\dots,y_n$ are variables such that $(y_1,\dots,y_n)$ is a permutation of the tuple $(x_1,\dots,x_n)$.

    \item The \emph{totally symmetric} condition $ \TS n $ contains all identities 
    $$
        \forall x_1,\dots,x_n\colon f(x_1,\dots,x_n)=f(y_1,\dots,y_n)
    $$
    where $x_1,\dots,x_n,y_1,\dots,y_n$ are variables with $\{x_1,\dots,x_n\}=\{y_1,\dots,y_n\}$.
    \item For every odd $n\ge 3$ the \emph{generalized minority} condition $\Gmin(n)$ contains the identity 
    $$
        \forall x_1,\dots,x_{n-2},y,z\colon f(x_1,\dots,x_{n-2},y,y)=f(x_1,\dots,x_{n-2},z,z)
    $$
    and contains additionally all identities in $\FS n$. 
\end{itemize}
\end{definition}

\subsection{Primitive Positive Constructions}\label{sec:ppconstructions}

 A \emph{primitive positive formula} (pp-formula) is a first order formula using only
 existential quantification and conjunctions of atomic formulas, i.e., formulas
 of the form $R(x_1,\dots,x_n)$ or $x=y$, and the special formula $\bot(x_1,\dots,x_n)$. The formula
 $\bot(x_1,\dots,x_n)$ always evaluates to the empty set. 
 Let $\A$ and $\B$ be relational structures. We say that $\B$ is a \emph{pp-power} of $\A$ if there is some $n\in \N^+$ such that $\B$ is
 isomorphic to a structure with domain $A^n$, whose relations
 are pp-definable in $\A$ (a k-ary relation on $A^n$ is regarded as a $kn$-ary relation
 on $A$).
 We say that $\A$ \emph{pp-constructs} $\B$, denoted by $\A\ppleq \B$, if $\B$ is
 homomorphically equivalent to a pp-power of $\A$. The structures $\A$ and $\B$ are \emph{pp-interconstructable}, denoted $\A\ppeq\B$, if $\A\ppleq\B$ and $\B\ppleq\A$. 
 It is well known that every finite structure $\structA$ is pp-interconstructable with a finite structure whose polymorphism clone is idempotent~\cite[Theorems 16–17]{Pol}. 
 
 The pp-constructability relation is a quasi order on the class
 of all relational structures~\cite{wonderland}. This quasi order can be characterized in different ways.

\begin{theorem}[{Theorem 1.3 in \cite{wonderland}}]\label{thm:ppconIffMinionHom}
    Let $\structA$ and $\structB$ be two finite structures. The following are equivalent:
    \begin{enumerate}
        \item The structure $\structA$ primitively positively constructs the structure $\structB$.
        \item There exists a minor preserving map $\Pol(\structA) \to \Pol(\structB)$.
        \item Every minor condition $\Sigma$ that is satisfied by $\pol(\A)$ is also satisfied by $\pol(\B)$.
    \end{enumerate}
\end{theorem}


\subsubsection{The Submaximal Element}
It is well known that in the poset of all finite structures ordered by pp-constructability $\C_1$ is a representative of the unique top element and $\P_1$ is a representative of its unique lower cover, see \cite[Proposition 2.11]{vucaj2024submaximal} or \cite[Proposition 3.2.5]{AlbertThesis}.
In~\cite{vucaj2024submaximal} Vucaj and Zhuk study moreover the subposet of all finite structures with at most three elements. They classify all lower covers of $\P_1$ in this poset. 
Some parts of their paper have been republished in the thesis~\cite{AlbertThesis} which was in the end published earlier.
To access the lower covers of $\P_1$ in the poset of all finite structures, we will use two results from the paper that were not limited to three element structures. The first result implies that generalized minority conditions of all odd arities and the totally symmetric conditions of all arities together form a generating set for the set of minor conditions satisfied by $\Pol(\P_1)$.

\begin{theorem}[Theorem 3.12 in~\cite{vucaj2024submaximal} and Theorem 6.1.13 in~\cite{AlbertThesis}]\label{Theorem:QuelleIdempotent}
    Let $\B$ be a finite structure. The following are equivalent:
    \begin{enumerate}
        \item The structure $\B$ is pp-interconstructable with $\bP_1$ or $\C_1$.
        \item The structure $\bP_1$ pp-constructs $\B$.
        \item The structure $\B$ has generalized minority polymorphisms of all odd arities and totally symmetric polymorphisms of all arities. 
    \end{enumerate}
\end{theorem}

The second result gives a sufficient condition for the existence of a quasi majority operation.

\begin{lemma}[{Lemma 3.4 in \cite{vucaj2024submaximal} and Lemma 6.1.5 in \cite{AlbertThesis}}] 
    \label{lem:SigmaPAndMaltImplyMajority}
    Let $\clone$ be a clone on a finite set with $k$ elements such that
         $\clone \models \Sigma_p$ for every prime $p \le k$ and
         $\clone \models \Malt$.
Then $\clone\models\Majority$.
\end{lemma}


%% file: 3_submaximal.tex
\section{The minor conditions satisfied by the clone of all idempotent operations}

To be able to access the lower covers of $\P_1$, we want to show a theorem about $\P_1$ first. 
The goal of this section is to show the following strengthening of Theorem~\ref{Theorem:QuelleIdempotent}.
\begin{theorem}\label{theorem:idempotence}\label{thm:idempotence}
    Let $\B$ be a finite structure such that $\B$ has a quasi Maltsev polymorphism and fully symmetric polymorphisms of all arities. 
    Then, $\B$ is pp-constructable from $\P_1$. 
\end{theorem}

This result is even stronger than a conjecture by Vucaj and Zhuk, which suggests that any finite structure with a quasi Maltsev polymorphism and totally symmetric polymorphisms can be pp-constructed from $\P_1$ (Conjecture~4.4 in~\cite{vucaj2024submaximal}).  

Combining Theorem~\ref{Theorem:QuelleIdempotent} with Lemma \ref{lem:SigmaPAndMaltImplyMajority}, it is clear that in order to prove Theorem~\ref{thm:idempotence} it suffices to show that any clone with a quasi majority, a quasi Maltsev, and fully symmetric operations of all arities also has generalized minority operations of all arities. 
As a first step we introduce two new families of minor conditions that will be important in intermediate steps in the proof.

\begin{definition} We consider some more minor conditions.
\begin{itemize}
    \item Let $n$ be odd. Then the \emph{compatible generalized minority} condition  in $(\gmin_1, \gmin_3,\gmin_5,\dots,\gmin_n)$ 
    contains for each odd $m$ in $[n-2]$ the identity 
    $$
        \forall x_1,\dots,x_{m},y\colon\gmin_{m+2}(y,y,x_1,\dots,x_m)=\gmin_m(x_1,\dots,x_m).
    $$
    and additionally it contains all identities of the fully symmetric condition in $\gmin_m$ for each odd $m$ in $[n]$.
    \item The \emph{\pairing{}} condition for $n\geq k\geq1$ with $n$ odd, denoted $\GPinter n k$, consists of the identities
    {
    \begin{align*}
        \forall x_1,\dots,x_{(n+1)/2}\colon f(y_i,y_i,\dots,y_i)= f(y_1,y_2,\dots,y_n) 
    \end{align*}
    for all $1\leq i\leq k$ 
    and all $y_1,\dots,y_n\in\{x_1,\dots , x_{(n+1)/2}\} $ such that each $x_m$ with $x_m\neq y_i$ occurs an even number of times in $(y_1,\dots,y_n)$.}


\item We write $\SymGPinter{n}$ for the condition that there is a fully symmetric $f$ with $f \models \GPinter{n}{n}$.
\end{itemize}
\end{definition}

Observe that a generalized minority operation  ignores any argument that occurs an even number of times and a generalized pairing operation only ignores arguments that occur an even number of times if there is only one argument which occurs an odd number of times and this argument occurs within the first $k$ entries.
Therefore, every generalized minority operation is also a generalized pairing operation. 
Note that a standard König's tree lemma argument shows that if a finite structure $\B$ satisfies the $n$-ary compatible generalized minority condition for all $n$, then there is an infinite chain $\gmin_1, \gmin_3, \gmin_5,\dots$ of polymorphisms of $\B$ such that for all odd $n$ we have that $\gmin_1,\gmin_3,\dots,\gmin_n$ satisfy the $n$-ary compatible generalized minority condition.
The next two lemmata are preparation for Lemma~\ref{Lemma:MajoMaltPairing}, which states that any clone with a quasi majority and a quasi Maltsev operation also has \pairing{} operations of all odd arities.





   \begin{lemma}\label{ClaimPairingA}
        Let $\clone$ be a clone, $n\in\N$ odd, and  $2\le k<n$. If $\clone\models\Majority$ and $\clone\models\GPinter n k$.  Then, $\clone\models\GPinter n {k+1}$.
    \end{lemma}
    \begin{proof}
    Let $\maj,\gp nk\in\clone$ such that $\maj\models\Majority$ and $\gp nk\models\GPinter n k$.
        It is easy to verify that the map that sends $(c_1,\dots,c_n)$ to
        $$
            \Majocolumns{
            \gp n {k}(c_1,\dots,c_{k-2},c_{k-1},c_{k\phantom{{}-1}},c_{k+1},c_{k+2},\dots ,c_n), \\
            \gp n {k}(c_1,\dots,c_{k-2},c_{k-1},c_{k+1},c_{k\phantom{{}-1}},c_{k+2},\dots ,c_n), \\
            \gp n {k}(c_1,\dots,c_{k-2},c_{k\phantom{{}-1}},c_{k+1},c_{k-1},c_{k+2},\dots ,c_n)
            }
        $$

        is an element of $\clone$ and satisfies $\GPinter n {k+1}$.
    \end{proof}

    \begin{lemma} \label{ClaimPairingB}
         Let $\clone$ be a clone and $n\geq3$ odd. If $\clone\models\Maltsev$ and $\clone\models\GPinter n n$ then $\clone\models\GPinter {n+2} {2}$.
    \end{lemma}
    \begin{proof}
    Let $\malt,p_n\in\clone$ such that $\malt\models\Maltsev$ and $p_n\models\GPinter nn$.
        It is easy to verify that the map that sends $(c_1,\dots,c_{n+2})$ to
        $$
            \Malcolumns{p_n(c_1,c_1,\dots,c_{1\phantom{{}-1}}),\\
            p_n(c_3,c_4,\dots,c_{n+2}) ,\\
            p_n(c_2,c_2,\dots,c_{2\phantom{{}-1}})
            }
        $$
        is in $\clone$ and satisfies $\GPinter {n+2}2$.
    \end{proof}

Combining Lemmata~\ref{ClaimPairingA} and~\ref{ClaimPairingB} in an induction argument we obtain the following lemma. Note that the Maltsev condition is equivalent to $\GPinter{3}{2}$.
\begin{lemma} \label{Lemma:MajoMaltPairing}
    Let $\clone$ be a clone with a quasi majority and quasi Maltsev operation. Then $\clone$ also has \pairing{} operations of every odd arity.
\end{lemma}



The next three lemmata show that any clone with generalized pairing operations of all odd arities and fully symmetric operations of all arities also has generalized minority operations of all odd arities and totally symmetric operations of all arities.

\begin{lemma} \label{lemma:gpsymmetrizable}
    The \pairing{} condition $\GPinter n n$ is symmetrizable, i.e., for every clone $\clone$ with $\clone\models \GPinter n n $ and $\clone\models\FS{n!}$ we have that $\clone\models \SymGPinter n$.
\end{lemma}
\begin{proof}
Let $f,p_n\in\clone$ with $p_n\models \GPinter n n$ and $f\models\FS{n!}$. It is easy to see that the function 
\begin{align*}
    (c_1,\dots,c_n)\mapsto f(p_n(c_{\sigma(1)},\dots,c_{\sigma(n)})\mid \sigma\in S_n)
\end{align*}
is in $\clone$ and satisfies $\SymGPinter n $.
\end{proof}

\begin{lemma} \label{Lemma:pairingGeneralizedmin}
    Let $\clone$ be an idempotent clone and $k\in\N^+$. If $\clone\models\SymGPinter {n}$ for all odd $n$, then $\mathcal C$ also has compatible generalized minorities $\gmin_1,\gmin_3,\dots,\gmin_n$ for all odd  $n$.
\end{lemma}
\begin{proof}
Let $p_1,p_2,\dots\in\clone$ with $p_1\models\SymGPinter{4^1-1},p_2\models\SymGPinter{4^2-1},\dots$. 
    Note that that the set $[2k+1]$ has exactly $4^k-1$ proper subsets of odd size.
    For $k\in \IN$, we define $\gmin_{2k+1}$ inductively: $\gmin_1(c)\coloneqq c$ and $\gmin_{2k+1}(c_1,\dots,c_{2k+1})$ is
    $$
    p_k(\gmin_{|A|}(c_a\mid a\in A) \mid A\subsetneq [2k+1] \text{ of odd size}).
    $$
    These maps are clearly fully symmetric. Assume that $\gmin_1,\dots,\gmin_{2k-1}$ are compatible symmetric generalized minorities. Let $c_1=c_2$. For any set $A\subsetneq[2k+1]$ of odd size let $A\setdifference\{1,2\}$ denote the symmetric difference of $A$ with $\{1,2\}$. Note that $|A\setdifference\{1,2\}|$ is also odd. If $A\setdifference\{1,2\}\neq[2k+1]$, then by assumption  
    \begin{align*}
        \gmin_{|A|} (c_a\mid a\in A) &= \gmin_{|A\setdifference\{1,2\}|}(c_a\mid a\in A\setdifference\{1,2\}).
    \end{align*}
    
    Hence, the only argument of $p_k$ that occurs an odd number of times is $\gmin_{2k-1}(c_3,\dots,c_{2k+1})$. Therefore, 
    \[\gmin_{2k+1}(c_1,c_1,c_3,\dots,c_{2k+1})=\gmin_{2k-1}(c_3,\dots,c_{2k+1})\]
    and $\gmin_1,\dots,\gmin_{2k-1},\gmin_{2k+1}$ are compatible generalized minorities. 
\end{proof}

\begin{lemma}\label{Lemma:generalizedminTotallysym}
    Let $n$ be odd and let $\clone$ be a clone with compatible generalized minority operations $\gmin_1,\gmin_3,\dots,\gmin_n$ and a fully symmetric operation of arity $2^{n-1}$. Then $\clone$ also has a totally symmetric operation of arity $n$.
\end{lemma}
\begin{proof}
    Let $f\in\clone$ with $f\models\FS{2^{n-1}}$. 
    Note that the set $[n]$ has exactly $2^{n-1}$ subsets of odd size. 
    We define
    $$
    t(c_1,\dots,c_{n})\coloneqq f(\gmin_{|A|}(c_a\mid a\in A) \mid A\subseteq [n] \text{ of odd size})
    $$
    Since $f$ and $\gmin_1,\dots,\gmin_n$ are fully symmetric so is $t$. 
    To show that $\ts n$ is totally symmetric define formally $c_2\coloneqq c_1$ and $(c'_1,c'_2,c'_3,\dots,c'_n)\coloneqq(c_1,c_3,c_3,\dots,c_n)$. Now, it suffices to prove
    \[\ts{n}(c_1,c_2,c_3,\dots,c_{n})=\ts{n}(c'_1,c'_2,c'_3,\dots,c'_{n})\]
    for all $c_1,c_3,\dots,c_n\in C$.
    Let $A\subseteq[n]$. If  $2\notin A$ then
    \[\gmin(c_a\mid a\in A)=\gmin(c'_a\mid a\in A)\]
    and if $2\in A$ then
    \begin{align*}
      \gmin(c_a\mid a\in A)
      &=\gmin(c_a\mid a\in A\setdifference\{1,2\})\tag{as $c_1=c_2$}\\
      &=\gmin(c'_a\mid a\in A\setdifference\{1,2\})\tag{as $2\notin A\setdifference\{1,2\}$}\\
      &=\gmin(c'_a\mid a\in A\setdifference\{1,2\}\setdifference\{2,3\})\tag{as $c'_2=c'_3$}\\
      &=\gmin(c'_a\mid a\in A\setdifference\{1,3\})  
    \end{align*}
    Since the map $A\mapsto A\setdifference\{1,3\}$ is a bijection from the odd subsets of $[n]$ that contain 2 onto itself and $f$ is fully symmetric we have that
    \[\ts{n}(c_1,c_1,c_3,\dots,c_{n})=\ts{n}(c_1,c_2,\dots,c_{n})=\ts{n}(c'_1,c'_2,\dots,c'_{n})=\ts{n}(c_1,c_3,c_3,\dots,c_{n})\]
    as desired. 
\end{proof}


Combining the previous lemmata we can prove the following theorem.
\begin{theorem}\label{theorem:idempotenceClonewise}
    Let $\clone$ be an idempotent clone with a quasi majority, a quasi Maltsev, and fully symmetric operations of all arities. Then $\clone$ also has compatible generalized minority operations of all odd arities and totally symmetric operations of all arities.
\end{theorem}
\begin{proof}
    By Lemma \ref{Lemma:MajoMaltPairing}, $\clone$ has \pairing{} operations. They are symmetrizable by Lemma~\ref{lemma:gpsymmetrizable}, so $\clone\models\SymGPinter n$ for all $n\in\N^+$. By Lemma~\ref{Lemma:pairingGeneralizedmin}, the clone also has compatible symmetric generalized minority operations of all odd arities and by Lemma~\ref{Lemma:generalizedminTotallysym} totally symmetric operations of all arities.
\end{proof}

Now we can proof the main theorem of this section.
\begin{proof}[Proof of Theorem~\ref{theorem:idempotence}]
    As the theorem considers a finite structure $\B$ we may assume that $\pol(B)$ is idempotent and can apply Lemma \ref{lem:SigmaPAndMaltImplyMajority} to get 
    $\Pol(\B)\models\Majority$. Now, we can apply Theorem~\ref{theorem:idempotenceClonewise} and get the existence of generalized minority polymorphisms and totally symmetric polymorphisms. By Theorem~\ref{Theorem:QuelleIdempotent}, we finally get that  $\B$ is pp-constructable from $\bP_1$.
\end{proof}

%% file: 4_Groups-or-cycles-paper.tex
\section{Finite simple group actions}
In this section we will give a description of the third layer of the poset of all finite structures ordered by pp-constructability. To do so we first introduce a new family of structures. A way to view any group action as a structure with binary signature.
\begin{definition}\label{def:S(GactX)}
    Let $G$ be a group acting on a set $X$. Define the structure $\S(G\curvearrowright X)$ that has domain $X$ and has for every $g\in G$ the binary relation 
    \[g^{\S(G\curvearrowright X)}=\{(x,gx)\mid x\in X\}.\]
    If the group is acting on itself we often abbreviate $\S(G\curvearrowright G)$ to $\S(G)$.
\end{definition}

Note that a $G$-structure homomorphism from $\S(G\curvearrowright X)$ to $\S(G\curvearrowright X')$ is exactly a $G$-set homomorphism from $X$ to $X'$.

\begin{lemma}
    Let $G$ be a group acting on a set $X$ and let $G$ be generated by $S$.  Then $\S(G\act X)$ is pp-interconstructable with the $S$-reduct of $\S(G\act X)$. 
\end{lemma}
\begin{proof}
    Let $\B$ be the $S$ reduct of $\S(G\act X)$. Clearly, $\S(G\act X)$ can pp-construct $\B$. For the other direction consider the set $P\subseteq \Sym(B)$ of all permutations on $B$ whose graphs can be pp-defined from $\B$. By definition this set contains the actions of all elements in $S$. Let $g_1,g_2\in G$ such that their actions on $B$ are contained in $P$ and let $\phi_1(x,y)$ and $\phi_2(x,y)$ by pp-definition of the actions of $g_1$ and $g_2$, respectively. Then $\exists z\colon \phi_1(z,y)\wedge\phi_2(z,x)$ is a pp-definition of the action of $g_1g_2^{-1}$. Clearly, $P$ also contains the action of the neutral element $e_G$. Hence, $P$ contains all actions in the subgroup generated by $S$ and $\B$ can pp-define $\S(G\act X)$.
\end{proof}
This lemma implies in particular that actions of cyclic groups are pp-interconstructable with directed graphs.

\begin{example}\label{exa:ZnStructure}
    Let $\Z/n\Z$ act on $\{0,\dots,n-1\}$ by $g(k)=g+k$. Then $\S(\Z/n\Z)$ is pp-interconstructable with $\C_n$. 

    Consider the cyclic group $G=\langle(1\,2)(3\,4\,5)\rangle$ acting on $\{1,\dots,5\}$. Then $\S(G\act\{1,\dots,5\})$ is pp-interconstructable with the disjoint union of $\C_2$ and $\C_3$. 
\end{example}

In this sense the actions of cyclic groups correspond exactly to disjoint unions of directed cycles which have been studied in \cite{smooth-digraphs}.

\begin{observation}
    Let $G$ be a cyclic group acting on a set $X$. Then $\S(G\act X)\ppeq \C$ for some disjoint union of directed cycles $\C$.
    Let $\C$ be a disjoint union of directed cycles, then there is a cyclic group $G$ acting on $C$ such that $\S(G\act C)\ppeq\C$.
\end{observation}

Note that the same group can act in a lot of different ways resulting in different structures. 
We now  introduce a special type of group action. Let $G$ be a finite group. A \emph{minimal fixed point free} action of $G$ on a set $X$ is an action that has no fixed point such that the action of each proper subgroup of $G$ on $X$ does have a fixed point.

\begin{definition}
    Let $G$ be a group and $\{X_i\mid i\in I\}$ the set of all primitive $G$-sets up to isomorphism. Define the $G$-set $\prim(G)$ as the set $\coprod_{i\in I} X_i$ with the simultaneous action of $G$.
\end{definition}
Note that, by Observation~\ref{obs:primitiveActionsOnMaximalSubgroups}, we can always assume without loss of generality that $\prim(G)$ consists of cosets of maximal subgroups of $G$ on which $G$ acts by left multiplication. Furthermore, if $G$ is finite, then $\prim(G)$ is also finite.
We now show that $G\act \prim(G)$ is a minimal fixed point free action
and that all structures associated to minimal fixed point free actions of $G$ are homomorphically equivalent.

\begin{lemma}\label{lem:minimalFixedPointFreeActionsAreUnionsOfPrimitiveActions}
    Let $G$ be a finite  group of order at least two acting on a set $X$. Then $\S(G\act X)$ is homomorphically equivalent to $\S(G\act \prim(G))$ if and only if the action of $G$ on $X$ is a minimal fixed point free action. 
    
    In particular, the action of $G$ on $\prim(G)$ is a minimal fixed point free action.
\end{lemma}
\begin{proof}
For the $\Rightarrow$-direction, assume that $\S(G\act X)$ and $\S(G\act \prim(G))$ are homomorphically equivalent. Note that $G\act \prim(G)$ has no fixed point as the point stabilizers of $G$ are maximal subgroups and not $G$ itself. As $\S(G\act X)$ has a homomorphism to $\S(G\act \prim(G))$, also $G\act X$ has no fixed point.


    Let $U$ be a proper subgroup of $G$. We want to show that the action of $U$ on $X$ has a fixed point. Hence, we assume without loss of generality that $U$ is a maximal proper subgroup. 
    Then $G/U$ is isomorphic (as $G$-set) to a component of $\prim(G)$. As for all $g\in U$ we have $gU=U$, this isomorphism defines an element $u\in \prim(G)$ such that $U$ acts trivial on $u$.
    Since $\prim(G)$ admits a $G$-set-homomorphism $f$ to $X$, there is an $x=f(u)\in X$ such that for every $g\in U$ the relation for $g$ in $\S(G\act X)$ has the loop $(x,x)$. Hence $x$ is a fixed point of the action of $U$ on $X$, as desired. In conclusion, $G\act X$  is a minimal fixed point free action.

    For the $\Leftarrow$-direction let $G\act X$ be a minimal fixed point free action. 
    We show the homomorphic equivalence componentwise. 
    Note that the connected components of $\S(G\act X)$ are exactly the $G$-orbits of $X$. 
    Let $x\in X$. 
    Since $G$ has no fixed points, the stabilizer of $x$ with respect to $G$ must be contained in some maximal proper subgroup $U_x$ of $G$. 
    The map $gx\mapsto gU_x$ is a well defined homomorphism from the connected component of $\S(G\act X)$ that contains $x$ to $\S(G\act\prim(G))$. 
    Hence, there is a homomorphism from $\S(G\act X)$ to $\S(G\act\prim(G))$. 
    Consider a connected component of $\S(G\act\prim(G))$. There is a maximal proper subgroup $U$ of $G$ such that this component is $G/U$, where $G$ acts by left multiplication. Since $G\act X$ is a minimal fixed point free action there is an $x_U\in X$ such that $U\act X$ has $x_U$ as a fixed point. The map $gU\mapsto gx_U$ is a well defined homomorphism from $G/U$ to $\S(G\act X)$. For well-definedness consider $gU=g'U$. Then $g^{-1}g'\in U$. Therefore, $gx_U=gg^{-1}g'x_U=g'x_U$, as desired.  
    Hence,  there also is a homomorphism from $\S(G\act\prim(G))$ to $\S(G\act X)$.
\end{proof}

For finite simple groups we can give an explicit description of $\prim(G)$.

\begin{lemma} \label{lem:descriptionPrimG}     
If $G$ is a finite abelian simple group, then $\prim(G)\iso G$ as $G$-sets.
    If $G$ is a finite, non-abelian simple group, then 
    \[\prim(G)\iso\{M< G\mid M\text{ maximal subgroup}\}\] 
    where the action is given by conjugation, i.e., $g.M=gMg^{-1}$ for $g\in G$.
    In both cases, the action is faithful.
\end{lemma}
\begin{proof}
    Every abelian finite simple group is isomorphic to $\IZ/p\IZ$ for some $p$ prime. Its only maximal subgroup is $\{0\}$ inducing its only primitive $\IZ/p\IZ$-set $\IZ/p\IZ$.

    For a non-abelian finite simple group $G$ consider a maximal subgroup $M$ and the group action of $G$ on the conjugacy class of $M$. The stabilizer $\stab(M)$ is $\{g\mid gMg^{-1}=M\}$ and satisfies $M\le \stab(M)\le G$. Since $M$ is not normal we have $M=\stab(M)\ne G$ and this action is primitive. It is isomorphic to the action of $G$ on $G/M$ by left-multiplication. Therefore, the conjugation action on all maximal subgroups contains all primitive actions. It has no duplicates, because the stabilizer of each point is a different maximal subgroup.

    The action is faithful, because the kernel $\ker(G\act\prim(G))$ is a proper normal subgroup of a simple group.
\end{proof}
In Section~\ref{sec:ExamplesOfnon-abelianFiniteSimpleGroupactions} we present  $\S(G\act\prim(G))$ for the three smallest non-abelian finite simple groups. 
In the rest of this section we will show that $\T_3$ together with the structures 
$$\{\S(G\act\prim(G))\mid G \text{ finite simple group}\}$$
are the lower covers of $\P_1$ in the poset of finite structures ordered by pp-constructibility. As a first step we show that structures associated to fixed point free actions of finite groups can always pp-construct a structure associated to a minimal fixed point free action of a finite simple group.


\begin{lemma} \label{lemma:groupsBelowSimple}
    Let $G$ be a finite group acting a finite set $X$. Then exactly one of the following is true
    \begin{enumerate}
        \item the action has a fixed point and $\S(G\act X)$ is pp-interconstructable with $\C_1$ or
        \item the action has no fixed point and there is a finite simple group $H$ such that $\S(G\act X)\ppleq\S(H\act \prim(H))$.
    \end{enumerate} 
\end{lemma}
\begin{proof}
    If $G$ has a fixed point $p$, then $X$ is homomorphically equivalent to the trivial $G$-set $G/G$ because the map that sends the single element of $G/G$ to $p$ is a homomorphism. Moreover, the loop $\S(G\act G/G)$ is clearly pp-interdefinable with the loop $\C_1$.
    
    Now we consider the case, that $G$ has no fixed point and is the smallest counter example (in the size of $G$). 
    Note that $\S(G\act X)$ can pp-construct $\S(U\act X)$ for any  subgroup $U$ of $G$. Hence, we can assume that every proper subgroup of $G$ has a fixed point. Otherwise just replace $G$ with this smaller subgroup.
    

    If $G$ is simple, then $\S(G\act X)$ is pp-interconstructable with $\S(G\act \prim(G))$ as they are homomorphically equivalent by Lemma~\ref{lem:minimalFixedPointFreeActionsAreUnionsOfPrimitiveActions}.
    If $G$ is not simple, then it has a proper nontrivial normal subgroup $N$. Let $X'\subseteq X$ be the set of points that are fixed by $N$. Since $N$ is a proper subgroup this set is non-empty. Since $N$ is a normal subgroup of $G$ we have that $X'$ is $G$-invariant, i.e., $G(X')=X'$. Hence, $X'$ is also a $G$-set by restriction of $G\act X$. Observe that $N$ acts trivially on $X'$. Therefore, $X'$ can also be seen as a $G/N$-set with the well defined action $(gN).x\coloneqq g.x$. 
    This action has no fixed point, as it would be a fixed point of the original $G$-action.
    The structure $\S(G \act X')$ is pp-constructable from $\S(G\act X)$ by shrinking the domain to $\{x\mid \bigwedge_{g\in N} g.x=x\}$. Since $\S(G\act X')$ and $\S(G/N\act X')$ have the same binary relations (just with different relation symbols) we have that $\S(G\act X)\ppleq\S(G/N\act X')$. The claim follows by induction on the size of $G$.
\end{proof}

\begin{lemma}\label{lem:noFSnImpliesCanPPconstSnAction}
    Let $\B$ be a finite relational structure with no fully symmetric polymorphism of arity $n$. 
    Then there is an action of $S_n$ on a finite set $X$ that has no fixed point such that $\B$ can pp-construct $\S(S_n\act X)$.
\end{lemma}
For more concise notation in the proof we will view functions from a set $I$ to $B$ as $I$-ary tuples in $B$ (instead of introducing an order on $I$ and considering $|I|$-ary tuples).  
\begin{proof}
    Let $\pol^{(n)}(\B)$ denote the set of $n$-ary polymorphisms of $\B$. Let $S_n$ act on this set by $\sigma.f=f_\sigma$. Since $\pol(\B)\not\models\FS n$ this action has no fixed point.  
    We will show that $\B$ can pp-construct $\S(S_n\act \pol^{(n)}(\B))$. 
     
    Consider the structure $\S$ with domain $B^{B^n}$ that has for each $\sigma\in S_n$ the binary relation 
    $\{(f,f_\sigma)\mid f\in\pol(\B)\}$. Clearly, $\S$ is homomorphically equivalent to $\S(S_n\act \pol^{(n)}(\B))$. 
    Note that $B^n$-tuples of $B$ are just maps from $B^n$ to $B$. To see why $\S$ is a $B^n$-th pp-power of $B$ let $\sigma\in S_n$. We define a formula $\phi_\sigma$ with variables $x_t$ and $y_t$ for every $t\in B^n$. Whenever $t_1,\dots,t_k$ are in $B^n$ such that their entries are componentwise in some relation $R^{\B}$, then add the conjunct $R(x_{t_1},\dots,x_{t_k})$ to $\phi_\sigma$. Furthermore, add for every $t\in B^n$ the conjunct $x_t=y_{t_\sigma}$. Note that any tuple in $\phi_\sigma^{\B}$ can be interpreted as two functions $f,g\colon B^n\to B$, where $f(t)$ is the value assigned to $x_t$ and $g(t)$ is the value assigned to $y_t$. Observe that the conjuncts of the form $R(x_{t_1},\dots,x_{t_k})$ ensure that $f$ is a polymorphism of $\B$ and the conjuncts of the form $x_t=y_{t_\sigma}$ ensure that $g=f_\sigma$. Therefore, $\phi_\sigma^{\B}$ is a pp-definition of $\{(f,f_\sigma)\mid f\in\pol^{(n)}(\B)\}$ in $\B$. Hence $\S$ is a $B^n$-th pp-power of $\B$, as desired.
\end{proof}

Combining Lemmata~\ref{lemma:groupsBelowSimple} and~\ref{lem:noFSnImpliesCanPPconstSnAction} we obtain the following theorem.

\begin{theorem}\label{lem:noFSnImpliesCanPPconstSimpleGroupaction}
    Let $\B$ be a finite relational structure with no fully symmetric polymorphism of arity $n$. Then there is a finite simple group $G$ such that $\B$ can pp-construct $\S(G\act\prim(G))$. 
    The group $G$ can be obtained from $S_n$ by repeatedly taking quotients or subgroups.
\end{theorem}

We also have the following well known theorem about $\T_3$. 

\begin{theorem}[Theorem 5.4 in \cite{starkeThesisDigraphsmoduloprimitivepositive}]
    Let $\B$ be a finite relational structure that has no quasi Maltsev polymorphism. Then $\B$ can pp-construct $\T_3$. 
\end{theorem}

Combining these two theorems with Theorem~\ref{thm:idempotence} we obtain the following corollary.
\begin{corollary}\label{cor:lowerCoversFirstVersion}
    Let $\B$ be a finite relational structure. If $\B<\P_1$, then $\B$ can pp-construct $\T_3$ or $\B$ can pp-construct $\S(G\act\prim(G))$ for some finite simple group $G$. 
\end{corollary}

\subsection{Incomparable group actions}
To show that all structures from Corollary~\ref{cor:lowerCoversFirstVersion} are lower covers of $\P_1$ we will show that they are pairwise incomparable with respect to pp-constructability. To do so we will associate a minor condition to any group action on a finite set. In order to make this association more natural we will, similar to clones, allow minor conditions to use $X$-ary function symbols for arbitrary finite sets $X$. 
\begin{definition}
    Let $X$ be a set and let $H\curvearrowright Y$ be a group action on a finite set $Y$. An operation $f\colon X^Y\to X$ satisfies the \emph{$H\curvearrowright Y$-condition}, denoted $\Sigma(H\curvearrowright Y)$, if 
        $f(t)= f(t_h)$ for all $h\in H$,  
    where $t_h\coloneqq Y\to X,y\mapsto t({h.y})$.
\end{definition}

\begin{example}\label{exa:ZnCondition}
    Let $\Z/n\Z$ act on $\{0,\dots,n-1\}$ by $g(k)=g+k$. Then an $n$-ary operation $f$ satisfies $\Sigma(\Z/n\Z)$ if and only if it satisfies $\Sigma_n$. Similarly, $f\models\Sigma(S_n\act[n])$ if and only if $f\models\FS n$. 
\end{example}

The following lemma characterises when a structure associated to an action satisfies the condition associated to an action.
\begin{lemma}\label{lem:GModlesHiffNoGlobalFixedPoint}
    Let $G$ and $H$ be finite groups acting on $X$ and $Y$, respectively. Then, the following are equivalent
    \begin{enumerate}
        \item $\Pol(\S(G\act X))\models \Sigma(H\act Y)$ and 
        \item for every $t\in X^Y$ there exists a point in $X$ that is a fixed point of $\stab_G(H(t))\le G$, where $H(t)=\{Y\to X,y\mapsto t({h.y}) \mid h\in H\}$ and the set stabilizer of $G$ is considered with respect to the componentwise action of $G$ on $X^Y$.
    \end{enumerate}
\end{lemma}

\begin{proof}
The lemma is proved by the same argument as \cite[Lemma 6.2]{smooth-digraphs}, which shows the statement for finite cyclic groups.

Assume that there is some $f\colon X^Y\to X$ in $\Pol(\S(G\act X))$ that satisfies $\Sigma(H\act Y)$ and consider $t\in X^Y$. As $f\models \Sigma(H\act Y)$, every point in the set $H(t)$ is mapped to $f(t)$. Thus, every  $g\in G$, that stabilizes $H(t)$, also stabilizes $f(t)$. Thus, the point $f(t)\in X$ is a fixed point of $\stab_G(H(t))$.

For the other direction assume that $\stab_G(H(t))$ has a fixed point for every $t\in X^Y$. We can now construct $f\colon X^Y\to X$ such that $f\models \Sigma(H\act Y)$.
First note that 
$$\{(t,t')\in X^Y\times X^Y\mid \exists g\in G, h\in H: t'=gt_h\}$$
is an equivalence relation on $X^Y$.
Now choose for every equivalence class a representative $t$ and choose $a_t$ as a fixed point of $\stab_G(H(t))$. Define furthermore $f(gt_h)$ as $ga_t$. By definition $f\models \Sigma(H\act Y)$. It is easy to verify that $f$ is a polymorphism of $\S(G\act X)$. To show that $f$ is well-defined, we need to show that if $gt_h=g't_{h'}$ are two representations of the same element, then
their images $ga_t$ and $g'a_t$ are equal.
We have 
\begin{align*}
&&g't_{h'}&=gt_h \\
\iff{}& & {g}^{-1}g't&=t_{h{h'}^{-1}} \\
\iff{}& \forall h''\in H\colon & {g}^{-1}g't_{h''}&=t_{h{h'}^{-1}h''} \in H(t) \\
\implies{}& & {g}^{-1}g' &\in \stab_G(H(t)) \\
\implies{}& &{g}^{-1}g'a_t&=a_t \\
\iff{}& &ga_t&=g'a_t
\end{align*}
and thus, $f$ is well-defined and $\Pol(\S(G\act X))\models \Sigma(H\act Y)$.  
\end{proof}

Using the previous lemma we can show that the structure associated to an action does not satisfy the condition associated to the same action.
\begin{lemma}\label{lem:GdoesnotsatisfySigmaG}
    Let $G$ be a finite group acting on a finite set $X$ without a fixed point. Then $\Pol(\S(G\act X))\not\models\Sigma(G\act X)$.
\end{lemma}
\begin{proof}
    Consider $t\coloneqq (X\to X,x\mapsto x)\in X^X$ and let $g\in G$. Then $(gt)(x)=g(t(x))=g(x)=t(g(x))=t_g(x)$ for all $x\in X$. Therefore, $g\in \stab_G(G(t))$. Hence, $\stab_G(G(t))=G$. By assumption $G\act X$ has no fixed point. Lemma~\ref{lem:GModlesHiffNoGlobalFixedPoint} implies that $\Pol(\S(G\act X))\not\models\Sigma(G\act X)$.
\end{proof}

We now show that the structure associated to a minimal fixed point free actions of a finite simple group satisfies the condition associated a minimal fixed point free actions of any other finite simple group. This shows that these structures cannot pp-construct each other.

\begin{lemma}\label{lem:SimpleGroupsModelSimpleConditions}
    Let $G$ and $H$ be finite simple groups. Then,
    \begin{align*}
        &\Pol(\S(G\curvearrowright\prim(G)))\models\Sigma(H\curvearrowright \prim(H)) &&\text{if and only if} &&G\not\iso H.
    \end{align*}
\end{lemma}
\begin{proof}
    If $G\iso H$, then, by Lemma~\ref{lem:GdoesnotsatisfySigmaG}, $\Pol(\S(G\curvearrowright\prim(G)))\not\models\Sigma(H\curvearrowright \prim(H))$.
    
    Conversely, if $\Pol(\S(G\curvearrowright\prim(G)))\not\models\Sigma(H\curvearrowright \prim(H))$, we want to use Lemma~\ref{lem:IsoOfSubquotient} to show $G\iso H$.
    By Lemma~\ref{lem:GModlesHiffNoGlobalFixedPoint}, there is a $t\in \prim(G)^{\prim(H)}$ such that $\stab_G(H(t))$ has no fixed point. Hence, by Lemma~\ref{lem:minimalFixedPointFreeActionsAreUnionsOfPrimitiveActions}, $\stab_G(H(t))=G$.
    By Lemma~\ref{lem:IsoOfSubquotient}, $\stab_G(t)$ is a normal subgroup of $\stab_G(H(t))=G$. As $G$ has no fixed point in $X$, $G$ also has no fixed point in $X^Y$ and thus $t$ is no fixed point. Therefore, $\stab_G(t)$ is a proper normal subgroup of the simple group $G$ and thus trivial.
    
    Assume that $\stab_H(G(t))$ would be a proper subgroup of $H$. Then, by Lemma~\ref{lem:minimalFixedPointFreeActionsAreUnionsOfPrimitiveActions}, it has a fixed point $y$ in $\prim(H)$. Let $g\in G$. Since $\stab_G(H(t))=G$, there is some $h\in H$ with $gt=t_h$. Hence, $gt(y)=t_h(y)=t(hy)=t(y)$. Therefore, $t(y)$ is a fixed point of $G$, a contradiction. Hence $\stab_H(G(t))=H$. 

    Thus, we get by Lemma~\ref{lem:IsoOfSubquotient} 
    $$
        G/\{1\} = \stab_G(H(t))/\stab_G(t)\iso \stab_H(G(t))/\stab_H(t) = H/\stab_H(t)
    $$
    and as $H$ is a simple group, the right side is $H$ or $\{1\}$. As $G\ne \{1\}$, we get $G\iso H$.
    \end{proof}

In order to finish our description of the lower covers of $\P_1$ we need to show that $\T_3$ and the structures associated to minimal fixed point free actions of finite simple groups cannot pp-construct each other. 
\begin{lemma}\label{lem:GroupModelsMaltsev}
    Let $G$ be a finite group acting on a finite set $X$. Then $\S(G\act X)$ has a quasi Maltsev polymorphism.
\end{lemma}
\begin{proof}
    Consider the minority operation $m\colon X^3\to X$ with $m(x,y,z)=x$ for all distinct $x,y,z\in X$. Clearly, $m\models\Malt$. Note that any relation of $\S(G\act X)$ is a disjoint union of cycles. Furthermore, every such relation is preserved by $m$. Hence, $m$ is a polymorphism of $\S(G\act X)$, as desired. 
\end{proof}

We need the following well known lemma.
\begin{lemma}\label{lem:T3modelsTSAndT3notModlesMaltsev}
    We have $\pol(\T_3)\not\models\Malt$. 
    Furthermore $\pol(\T_3)\models \TS n$ for all $n\in \IN$. In particular, $\pol(\T_3)\models \Sigma(G\act X)$ for all finite groups $G$ acting on a finite set $X$. 
\end{lemma}
\begin{proof}
    Let $m$ be a quasi Maltsev operation over $\{0,1,2\}$. Then 
    $m(1,0,0)=m(1,1,1)=m(2,2,1)$. Since $((1,0,0),(2,2,1))$ is an edge in $\T_3^3$ we have that $m$ can not be a polymorphism of $\T_3$.
    The $n$-ary minimum operation on $\{0,1,2\}$ is totally symmetric and a polymorphism of $\T_3$. Hence $\pol(\T_3)\models \TS n$.
\end{proof}



\begin{theorem}\label{thm:mainLowerCoversClassification}
    The structures $$\{\bT_3\} \cup \{\S(G\act\prim(G))\mid G \text{ finite simple group}\}$$ are representatives of the lower covers of $\bP_1$ in the poset of all finite structures ordered by pp-constructability. In particular, a finite structure $\B$ is not pp-constructable from $\bP_1$ if and only if $\B$ can pp-construct one of the lower covers of $\bP_1$.
\end{theorem}
\begin{proof}
By Corollary~\ref{cor:lowerCoversFirstVersion}, any structure strictly below $\P_1$ is below on of these structures.
And by Lemmata~\ref{lem:SimpleGroupsModelSimpleConditions}, \ref{lem:GroupModelsMaltsev}, and~\ref{lem:T3modelsTSAndT3notModlesMaltsev} the set of minor conditions
\[\{\Malt\}\cup \{\Sigma(G\act\prim(G))\mid G \text{ finite simple group}\}\]
shows that these structures are all pairwise incomparable. 
\end{proof}

Combining Theorem~\ref{thm:mainLowerCoversClassification} with Lemmata~\ref{lem:SimpleGroupsModelSimpleConditions} and~\ref{lem:GdoesnotsatisfySigmaG} we obtain the following strengthening of Theorem~\ref{thm:idempotence}. 
\begin{corollary}\label{cor:weakerConditionBasisForP1}
    Let $\B$ be a finite structure such that $\Pol(\B)$ has a quasi Maltsev polymorphism and polymorphisms satisfying $\Sigma(G\act\prim(G))$ for every finite simple group $G$. Then $\B$ is pp-constructable from $\P_1$. 
\end{corollary}



Note that it does not follow from our previous results that $\S(G\act\prim(G))$ is a blocker for $\Sigma(G\act\prim(G))$. We will now show this this is the case, even for groups $G$ that are not simple.
\begin{theorem}
    Let $\B$ be a finite structure and $G$ be a finite group. Then $\Pol(\B)\not\models\Sigma(G\act\prim(G))$ if and only if $\B$ pp-constructs $\S(G\act\prim(G))$.
\end{theorem}
\begin{proof}
    The only-if-direction follows from Lemma~\ref{lem:GdoesnotsatisfySigmaG} as pp-constructions preserve minor identities.

    For the if-direction, let $\Pol^{(\prim(G))}(\B)$ denote the $\prim(G)$-ary polymorphisms of $\B$. 
    Consider the structure $\S$ with domain  $B^{B^{\prim(G)}}$ that has for each $g\in G$ the binary relation $\{(f,g.f)\mid f\in\Pol^{(\prim(G))}(\B)\}$, where $g.f=f_g\colon B^{\prim(G)}\to B, t\mapsto f(t_g)$ is the minor. 
    The proof that $\B$ can pp-construct $\S$ is analogous to the proof of Lemma~\ref{lem:noFSnImpliesCanPPconstSnAction}. 
    We now show that $\S$ and $\S(G\act\prim(G))$ are homomorphically equivalent. 
    It is easy to see that the map $\prim(G)\to B^{B^{\prim(G)}}, x\mapsto \pi_x$ is a homomorphism, where $\pi_x\colon B^{\prim(G)}\to B, t\mapsto t_x$.
    
    For the other homomorphism consider a connected component $\C$ of $\S$ that is not an isolated point. The structure $\C$ encodes a transitive action of $G$ on $C$. Let $f\in C$ and let $U_f$ be the stabelizer of $f$ with respect to $G$. Since $\Pol(\B)\not\models\Sigma(G\act\prim(G))$ we have that $U_f$ is not $G$. Hence, $U_f$ is contained in a maximal proper subgroup $U$ of $G$. It is easy to verify that the map $g.f\mapsto gU$ is a homomorphism from $\C$ to $\S(G\act\prim(G))$. Repeating this for every connected component of $\S$ yields a homomorphism from $\S$ to $\S(G\act\prim(G))$ as desired. 
\end{proof}

Using this theorem we can show that our set of minor conditions is the only generating set for the set of minor conditions satisfied in $\Pol(\P_1)$, where no condition is decomposable into weaker conditions. 
\begin{corollary}\label{cor:weakestGeneratingSet}
    Let $\conditionset$ be a set of minor conditions such that for all finite structures $\structA$ we have that $\Pol(\structA) \models \Sigma$ for all $\Sigma\in \conditionset$ implies that $\P_1$ pp-constructs $\structA$. Then every condition in
    $$
        \{\Malt\} \cup \{\Sigma(G\act\prim(G))\mid G \text{ finite simple group}\}
    $$
    is implied (for clones over a finite domain) by a single condition in $\conditionset$.
\end{corollary}
\begin{proof}
    Consider any condition $\Sigma'$ in 
    \[\{\Malt\} \cup \{\Sigma(G\act\prim(G))\mid G \text{ finite simple group}\}\] 
    and let $\structA$ be the corresponding blocker structure $\T_3$ or $\S(G\act \prim(G))$.
    Then, as $\P_1$ cannot pp-construct $\structA$, we get that $\Pol(\structA) \not\models \Sigma$ for one $\Sigma\in \conditionset$. So we get for all finite structures $\structB$ the implications
    $$
        \Pol(\structB) \not\models \Sigma' \Leftrightarrow \structB \le \structA \Rightarrow \Pol(\structB) \not\models \Sigma
    $$
    hence $\Sigma$ implies $\Sigma'$ (for clones over a finite domain).
\end{proof}
Note that this proof works whenever we have a set of minor conditions such that these conditions have blocker structures and these blocker structures are the lower covers of some structure.

%

%% file: A_Simple-group-actions.tex
\section{Examples of Non-Abelian Simple Group Actions}\label{sec:ExamplesOfnon-abelianFiniteSimpleGroupactions}
In this section we will present the structures associated to minimal fixed point free group actions of the three smallest non-abelian finite simple groups. By~\cite{smallSimpleGroups}, these three groups are
\begin{itemize}
    \item $A_5$, the alternating group of degree five with 60 elements,
    \item $\operatorname{PSL}(2,7)$, the projective special linear group of degree two over the seven element field with 168 elements, and
    \item $A_6$, the alternating group of degree six with 360 elements.
\end{itemize}

\subsection{The Alternating Group of Degree Five}
The smallest non-abelian finite simple group is the alternating group of degree five, denoted $A_5$. It has 60 elements and is the subgroup of $S_5$ generated by  the two permutation $f=(3\,4\,5)$ and $g=(1\,3)(2\,4)$.  Our goal is to determine the structure $\S(A_5\act\prim(A_5))$. Since $A_5$ is a permutation group we already have an action of $A_5$ on the set $[5]$. The following lemma shows that this action is primitive. 
\begin{lemma}\label{lem:minimalTransitiveActionIsPrimitive}
    Let $G$ be a finite group and let $n\in\N\setminus\{0,1\}$ be minimal such that $G$ can act transitively on $[n]$. Then any transitive action of $G$ on $[n]$ is primitive.
\end{lemma}
\begin{proof}
    Let $P$ be a partition of $[n]$ that is respected by $G\act[n]$. Then $G$ also acts transitively on $P$. Since $n$ is minimal in $\N\setminus\{0,1\}$ we have that $|P|$ is either 1 or $n$. In the  first case $P=\{[n]\}$ and in the second case $P=\{\{k\}\mid k\in[n]\}$. Hence $G\act[n]$ is primitive.   
\end{proof}

\begin{figure}
    \centering
    \Twentyone
    \caption{The $\{f,g\}$-reduct of $\S(A_5\act\prim(A_5))$, where the relation for $g$ is represented with undirected dashed edges.}
    \label{fig:A5}
\end{figure}

Note that since $A_5$ is simple any transitive action of $A_5$ on a set $X$ with $|X|\geq2$ must give an embedding of $A_5$ into $\Sym(X)$. 
Hence, Lemma~\ref{lem:minimalTransitiveActionIsPrimitive} can be applied to the action of $A_5$ on $[5]$ since its clearly transitive and there is no $2\leq n<5$ such that $A_5$ embeds into $S_n$ as $|A_5|=60>24=|S_4|$. Therefore, the action of $A_5$ on $[5]$ is primitive and $\S(A_5\act[5])$ is isomorphic to a connected component of $\S(A_5\act\prim(A_5))$, see Figure~\ref{fig:A5} (left) for the $\{f,g\}$-reduct of this component. To determine the  other connected components we need the conjugacy classes of the maximal subgroups of $A_5$. We used a python program to determine that there are three conjugacy classes of maximal subgroups that can be represented by 
\begin{itemize}
    \item the 12 element subgroup $U_1=\langle(1\,3)(2\,5),(1\,3\,5)\rangle$, 
\item the 6 element subgroup $U_2=\langle(1\,3)(2\,4),(1\,5)(2\,4)\rangle$, and 
\item the 10 element subgroup $U_3=\langle(1\,4)(2\,5),(1\,3\,4\,2\,5)\rangle$. 
\end{itemize}
In Figure~\ref{fig:A5} we have the three connected components of the $\{f,g\}$-reduct of $\S(A_5\act\prim(A_5))$ as computed by our program. The action can be seen as an embedding of $A_5$ into $S_{21}$ that maps 
\begin{align*}
    f&\text{ to } (3\,4\,5)(7\,8\,9)(10\,11\,12)(13\,14\,15)(16\,17\,18)(19\,20\,21)\text{ and}\\
    g&\text{ to } (1\,3)(2\,4)(6\,7)(8\,15)(9\,10)(12\,13)(18\,19)(20\,21). 
\end{align*}
This corresponds to the action of $A_5$ on maximal subgroups of $A_5$ by conjugation, here $U_1, U_2$, and $U_3$ correspond to the vertices labelled 1, 6, and 16, respectively. 
Note that the action of $A_5$ on $A_5/U_i$ by left-multiplication is isomorphic to the action on the conjugacy classes of $U_i$ by conjugation. Hence, the connected component of $\S(A_5\act\prim(A_5))$ containing $U_i$ consists of $[A_5:U_i]$ many elements.

We are not the first to study the structure $\S(A_5\act\prim(A_5))$. In a paper from 2016 Carvalho and Krokhin give this structure as an example of a structure with $p$-cyclic polymorphisms of all arities but no fully symmetric polymorphism of arity 5~\cite{carvalhoKrokhin21ElementStructure}. Their paper was very important for us since it showed that the list $\T_3,\C_2,\C_3,\C_5,\dots$ does not contain all lower covers of $\P_1$. However, this 21-element structure was quite mysterious to us in beginning. As they did not draw the connection to finite simple groups (it was not of interest to them) and it was unclear whether there were other structures and whether the structure they found was itself a lower cover of $\P_1$ or just below some up to this point unknown lower cover of $\P_1$.

The main result about $\S(A_5\act\prim(A_5))$ in the paper of Carvalho and Krokhin is the following theorem.

\begin{theorem}[Theorem 3 in \cite{carvalhoKrokhin21ElementStructure}]
The structure $\S(A_5\act\prim(A_5))$ has cyclic polymorphisms of all arities, but no fully symmetric polymorphism of arity 5.
\end{theorem}
We reproved the first part of this theorem in more generality in Lemma~\ref{lem:SimpleGroupsModelSimpleConditions}. The second part is easy since any 5-ary fully symmetric polymorphism would have to map $(1,2,3,4,5)$ to an element with an $f$ and a $g$ loop, which does not exist. For the same reason $\S(A_5\act\prim(A_5))$ does not have a fully symmetric polymorphism of any arity that can be written as a sum of 5's, 10's and 6's, i.e., any arity in $\{5,6,10,11,12,15,16,17,18\}\cup\{n\in\N\mid n\geq 20\}$.
In general we make the following observation.

\begin{observation}\label{obs:FS}
    Let $G$ be a finite group acting on a finite set $X$. Then
$\Pol(\S(G\act X))\not\models \FS k$ for all $k$ that can be written as the sum of sizes of connected components of $\S(G\act X)$.
\end{observation}

For structures associated to a minimal fixed point free action the converse is also true.
\begin{lemma}
    Let $G$ be a finite group. Then $\Pol(\S(G\act \prim(G)))\models \FS k$ if and only if $k$ can not be written as the sum of sizes of connected components of $\S(G\act \prim(G))$.
\end{lemma}
\begin{proof} 
    The only if direction is Observation~\ref{obs:FS}. 
    For the if direction note that $\FS k$ is equivalent to $\Sigma(S_k\act[k])$. Hence it suffices to verify the condition from Lemma~\ref{lem:GModlesHiffNoGlobalFixedPoint}. 
    Let $t\in \prim(G)^{[k]}$. We have to show that $\stab_G(S_n(t))$ has a fixed point. Let $a,b\in \prim(G)$ be in the same connected component such that the sets $\{i\in[k]\mid t_i=a\}$ and $\{i\in[k]\mid t_i=b\}$ are of different sizes. Such $a$ and $b$ exist, because $k=|\{i\in [k]\}|$ is not the sum of sizes of connected components of $\prim(G)$. Let $g\in G$ with $ga=b$. By choice of $a$ and $b$ we obtain that $gt$ and $t$ have a different number of entries equal to $b$. Therefore, $gt$ is not in the orbit $S_k(t)$ and $g\notin\stab_G(S_n(t))$. Hence, $\stab_G(S_n(t))$ is a proper subgroup of $G$. By Lemma~\ref{lem:minimalFixedPointFreeActionsAreUnionsOfPrimitiveActions}, $\stab_G(S_n(t))$ has a fixed point.
    %
\end{proof}

Combining this Lemma with Lemma~\ref{lem:minimalTransitiveActionIsPrimitive} yields the following corollary.
\begin{corollary}
        Let $G$ be a finite simple group. Then, the following numbers coincide:
    \begin{enumerate}
        \item the smallest $n$ such that $\Pol(\S(G\act\prim(G)))\not\models \FS n$,
        \item the size of a smallest connected component of $\S(G\act\prim(G))$, 
        \item the largest index of a maximal proper subgroup of $G$, and
        \item the smallest $n$ such that $G$ has an embedding into $S_n$.
    \end{enumerate}
\end{corollary}






\subsection{The Projective Special Linear Group of Degree Two Over The Seven Element Field}
The second smallest non-abelian group is $\operatorname{PSL}(2,7)$. It has 168 elements and is isomorphic to the subgroup of $S_7$ generated by the two permutations $f=(1\,2\,3\,4\,5\,6\,7)$ and $g=(2\,6)(3\,4)$. 

Our program computed that there are three conjugacy classes of maximal subgroups that can be represented by 
\begin{itemize}
    \item the 24 element subgroup $U_1=\langle(2\,7\,6\,5)(3\,4),(2\,7\,3)(4\,6\,5)\rangle$, 
\item the 21 element subgroup $U_2=\langle(1\,2\,4)(3\,6\,5),(1\,7\,6\,5\,4\,3\,2)\rangle$, and 
\item the 24 element subgroup $U_3=\langle(1\,5)(4\,6),(1\,7)(2\,3\,4\,6)\rangle$. 
\end{itemize}

In Figure~\ref{fig:PSL27} we have the three connected components of the $\{f,g\}$-reduct of $\S(\operatorname{PSL}(2,7)\act\prim(\operatorname{PSL}(2,7)))$ as computed by our program. The action can be seen as an embedding of $\operatorname{PSL}(2,7)$ into $S_{22}$ that maps 
\begin{align*}
    f&\text{ to } (1\,2\,3\,4\,5\,6\,7)(9\,10\,11\,12\,13\,14\,15)(16\,17\,18\,19\,20\,21\,22)\text{ and}\\
    g&\text{ to } (2\,6)(3\,4)(8\,9)(10\,15)(11\,12)(13\,14)(17\,21)(19\,18). 
\end{align*}
This corresponds to the action of $\operatorname{PSL}(2,7)$ on maximal subgroups of $\operatorname{PSL}(2,7)$ by conjugation, here $U_1, U_2$, and $U_3$ correspond to the vertices labelled 1, 8, and 22, respectively.


\begin{figure}
    \centering
\begin{tikzpicture}
    \node at (-4.25,0) {\TwentytwoB};
    \node at (0,0) {\TwentytwoA};
    \node at (4.25,0) {\TwentytwoC};
\end{tikzpicture}

    \caption{The $\{f,g\}$-reduct of $\S(\operatorname{PSL}(2,7)\act\prim(\operatorname{PSL}(2,7)))$, where the relation for $g$ is represented with undirected dashed edges.}
    \label{fig:PSL27}
\end{figure}

\subsection{The Alternating Group of Degree Six}
The third smallest non-abelian group is $A_6$. It has 360 elements and can be generated by the two permutations 
$f=(1\,2\,3\,5)(4\,6)$ and 
$g=(1\,2)(3\,4)$. 


Our program computed that there are five conjugacy classes of maximal subgroups that can be represented by 
\begin{itemize}
    \item the 60 element subgroup $U_1=\langle(2\,3)(4\,5),(1\,2\,4\,3\,5)\rangle$, 
\item the 36 element subgroup $U_2=\langle(1\,5\,3\,2)(4\,6),(1\,3\,6)(2\,5\,4)\rangle$,  
\item the 60 element subgroup $U_3=\langle(1\,3\,5)(2\,4\,6),(1\,2\,3\,6\,4)\rangle$,
\item the 24 element subgroup $U_4=\langle(1\,2)(3\,4\,5\,6),(1\,2)(3\,6\,4\,5)\rangle$, and  
\item the 24 element subgroup $U_5=\langle(1\,4)(2\,3),(1\,2\,5)(3\,4\,6)\rangle$.
\end{itemize}

In Figure~\ref{fig:A6} we have the five connected components of the $\{f,g\}$-reduct of $\S(A_5\act\prim(A_5))$ as computed by our program. The action can be seen as an embedding of $A_6$ into $S_{52}$ that maps 
\begin{align*}
    f\text{ to } &(1\,2\,3\,5)(4\,6)(8\,9\,10\,11)(12\,13\,14\,15)(17\,18)(19\,20\,21\,22)(23\,24\,25\,26)\\
    &(27\,28\,29\,30)(31\,32\,33\,34)(35\,36)(38\,39\,40\,41)(42\,43\,44\,45)(46\,47\,48\,49)(50\,51)\\
    &\text{ and}\\
    g\text{ to } &(1\,2)(3\,4)(7\,8)(10\,12)(11\,13)(15\,16)(18\,19)(20\,21)(23\,34)(24\,35)(26\,27)(28\,31)\\
    &(32\,36)(33\,37)(38\,49)(39\,52)(40\,50)(41\,42)(43\,46)(48\,51). 
\end{align*}
This corresponds to the action of $A_5$ on maximal subgroups of $A_5$ by conjugation, here $U_1, U_2$, $U_3$, $U_4$, and $U_5$ correspond to the vertices labelled 6, 7, 17, 25 and 47, respectively.


\begin{figure}
    \centering

\FiftytwoC \FiftytwoA \FiftytwoCMirrow

\FiftytwoB~\FiftytwoBMirrow

    \caption{The $\{f,g\}$-reduct of $\S(A_6\act\prim(A_6))$, where the relation for $g$ is represented with undirected dashed edges.}
    \label{fig:A6}
\end{figure}



\subsection{Final Observations}


Note that in the $\{f,g\}$-reducts of all three examples we have that for each connected component $\C$ there is a \emph{dual} component, i.e., a component that can be obtained from $\C$ by reversing all edges. Often the dual component is isomorphic to the original one. That  every component has a dual component comes from the fact that  in these three examples the map $f\mapsto f^{-1}, g\mapsto g^{-1}$ can be extended to an automorphism of $G$. 
However, this is not true in general.
Thanks to Dave Benson for providing the following counterexample. 
Consider the sporadic group $M_{11}$, generated by an element $f$ of order 8 and an element $g$ of order 11. By~\cite{Wilson1985ATLASOF}, no maximal subgroup can contain these elements, so they do indeed generate $M_{11}$. There are no outer automorphisms of $M_{11}$, hence all automorphisms of $M_{11}$ are inner, i.e., they are conjugation by some fixed element. Since in $M_{11}$ an element of order 8 is not conjugate to its inverse, there is no automorphism sending $f$ to its inverse. However, it could be that there is a different choice of generators that works. But this is not always the case, see~\cite{mathoverflowGroupAutomophism} for an example of a group where every set of two generators fails to yield an automorphism.

So far every lower cover of $\P_1$ can be represented by a digraph or a structure with two binary relations. This is true in general as it is well known that every finite simple group can be generated by two elements. The proof of this statement depends on the classification of finite simple groups.




\subsection*{Acknowledgements}
The authors want to thank Manuel Bodirsky and Andrew Moorhead for thoroughly reading our draft and for providing man helpful comments. 